\documentclass[times]{fldauth}

\usepackage{moreverb}
\usepackage{graphicx}
\usepackage{amssymb}
\usepackage{amsmath}
\usepackage{nicefrac} 
\usepackage{dsfont}

\usepackage[colorlinks,bookmarksopen,bookmarksnumbered,citecolor=red,urlcolor=red]{hyperref}

\newcommand\BibTeX{{\rmfamily B\kern-.05em \textsc{i\kern-.025em b}\kern-.08em
T\kern-.1667em\lower.7ex\hbox{E}\kern-.125emX}}

\usepackage{color} 
\usepackage{psfrag}

\newcommand{\halb}{{\frac{1}{2}}}
\newcommand{\CFL}{{\textnormal{CFL}}}
\newcommand{\Q}{{\mathbf{Q}}}
\newcommand{\f}{{\mathbf{f}}}

\def\volumeyear{2018}

\begin{document}

\runningheads{M. Dumbser et al.}{Semi-implicit divergence-free finite volume schemes for MHD}

\title{A divergence-free semi-implicit finite volume scheme for ideal, viscous and resistive magnetohydrodynamics}

\author{M. Dumbser$^1$\corrauth, D.S. Balsara$^2$, M. Tavelli$^1$, F. Fambri$^1$}

\address{\center $^1$ Department of Civil and Environmental Engineering, University of Trento, Via Mesiano, 77, 38123 Trento, Italy \\   
          $^2$ Physics Department, University of Notre Dame du Lac, 225 Nieuwland Science Hall, Notre Dame, IN 46556, USA}

\corraddr{michael.dumbser@unitn.it, dbalsara@nd.edu, m.tavelli@unitn.it, francesco.fambri@unitn.it}

\begin{abstract}
In this paper we present a novel \textit{pressure-based} semi-implicit finite volume solver for the equations of \textit{compressible} ideal, viscous and resistive magnetohydrodynamics 
(MHD). The new method is conservative for mass, momentum and total energy and in multiple space dimensions it is constructed 
in such a way as to respect the divergence-free condition of the magnetic field exactly, also in the presence of resistive effects. This is possible via the use
of multi-dimensional Riemann solvers on an appropriately staggered grid for the time evolution of the magnetic field and a double curl formulation of the resistive terms. 
The new semi-implicit method for the MHD equations proposed here discretizes the nonlinear convective terms as well as the time evolution of the magnetic field  
\textit{explicitly}, while all terms related to the pressure in the momentum equation and the total energy equation are discretized \textit{implicitly}, making again use 
of a properly staggered grid for pressure and velocity. 
Inserting the discrete momentum equation into the discrete energy equation then yields a \textit{mildly nonlinear} symmetric and positive definite algebraic 
system for the pressure as the only unknown, which can be efficiently solved with the (nested) Newton method of Casulli et al. The pressure system becomes 
\textit{linear} when the specific internal energy is a linear function of the pressure. 
The time step of the scheme is restricted by a CFL condition based only on the fluid velocity and the Alfv\'en wave speed and is not based on the speed of the magnetosonic 
waves. 
Being a semi-implicit pressure-based scheme, our new method is therefore particularly well-suited for low Mach number flows and for the incompressible limit of the  
MHD equations, for which it is well-known that explicit \textit{density-based} Godunov-type finite volume solvers become increasingly inefficient and 
inaccurate due to the increasingly stringent CFL condition and the wrong scaling of the numerical viscosity in the incompressible limit. 
We show a relevant MHD test problem in the low Mach number regime where the new semi-implicit algorithm is a factor of 50 faster than a traditional explicit  
finite volume method, which is a very significant gain in terms of computational efficiency. However, our numerical results confirm that our new method 
performs well also for classical MHD test cases with strong shocks. In this sense our new scheme is a true \textit{all Mach number} flow solver.   
\end{abstract}


\maketitle

\section{Introduction}
\label{sec.intro}

Since their invention by Harlow and Welch \cite{markerandcell}, pressure-based semi-implicit finite difference schemes on staggered grids have become widespread 
over the last decades for the solution of the incompressible Navier-Stokes equations with and without moving free surface, see e.g.  
\cite{chorin1,chorin2,patankar,patankarspalding,BellColellaGlaz,vanKan,HirtNichols,Casulli1990,CasulliCheng1992,Casulli1999,CasulliWalters2000,Casulli2009,CasulliVOF}
for a non-exhaustive overview of some of the most important contributions. An early application of semi-implicit schemes to compressible gas dynamics was the method of 
Casulli and Greenspan \cite{CasulliCompressible}, but their scheme was not conservative and therefore unable to solve problems including shock waves. 
In the field of numerical methods for high Mach number compressible flows, typically explicit density-based Godunov-type finite volume schemes 
\cite{lax,godunov,roe,osherandsolomon,hll,munz91,munz94,Toro:1994,LeVeque:2002a,toro-book} are preferred, since they are by construction conservative and thus allow  
the correct computation of shock waves. Therefore, the application of semi-implicit methods to compressible flows with shock waves is still quite rare, and some recent 
developments in this direction have been made only very recently in \cite{MunzPark,CordierDegond,DumbserCasulli2016,RussoAllMach}, where new \textit{conservative} 
pressure-based semi-implicit schemes have been proposed that are also suitable for the simulation of flow problems including shock waves. 
Concerning the numerical simulation of compressible magnetized plasma flows governed by the ideal or viscous and resistive magnetohydrodynamics (MHD) equations, 
only very little work has been done so far concerning the development of 
semi-implicit schemes. The existing semi-implicit schemes for MHD either apply only to the incompressible case, or they are not based on a conservative formulation, 
see e.g. \cite{Amari,Lerbinger,Harned,Finan}. The declared aim of this paper is therefore to close this gap and to propose a new conservative and pressure-based 
semi-implicit finite volume method for the solution of the compressible ideal and viscous and resistive MHD equations that applies both to high Mach number flows 
with shocks as well as to low Mach number or even incompressible flows. It is well-known that explicit density-based solvers become increasingly inefficient and 
inaccurate in the low Mach number regime and therefore an implicit time discretization is needed. However, discretizing all terms implicitly would lead to a 
\textit{highly nonlinear} non-symmetric system with a large number of unknowns (density, velocity, pressure and magnetic field), for which convergence is very 
difficult to control. Therefore, the new semi-implicit finite volume (SIFV) method proposed in this paper uses instead an explicit discretization for all nonlinear 
convective terms and for the time evolution of the magnetic field, while an implicit discretization is only employed for the pressure terms. This judicious combination 
leads in the end to only one \textit{mildly-nonlinear} and \textit{symmetric positive definite} system for the fluid pressure as the only unknown. 
The properties of the pressure system allow the use of the Newton-type techniques of Casulli et al. 
\cite{BrugnanoCasulli,BrugnanoCasulli2,CasulliZanolli2010,CasulliZanolli2012}, for which convergence has been \textit{rigorously proven}. 
Due to the implicit pressure terms, the time step of our new scheme is only restricted by the fluid velocity and the Alfv\'en wave speed, and not by the speeds 
of the magnetosonic waves. For this reason, the method proposed in this paper is a true \textit{all Mach number} flow solver. 
 
Modern computer codes for the solution of the MHD equations are mainly based on second or higher order Godunov-type finite volume schemes 
\cite{fallemhd3,Balsara2004,JiangWu,GardinerStone,PowellMHD2,Dumbser2008,AMR3DCL,LagrangeMHD,Zanotti2015c,ADERdivB} or on the discontinuous Galerkin (DG) finite element framework 
\cite{WarburtonVRMHD,Zanotti2015c,divdg1,divdg2,divdg3,divdg4,EntropyDGMHD}. In all these methods, the proper discretization of the magnetic field is of fundamental 
importance due to the well-known divergence-free constraint which the magnetic field must satisfy. Various solutions to this problem have been proposed 
in the literature so far and they can be essentially classified in two main categories: 
i) the first class contains the exactly divergence-free methods, following the ideas of Balsara and Spicer \cite{BalsaraSpicer1999,Balsara2004} and which requires the 
electric field at the vertices of each element and thus a multi-dimensional Riemann solver \cite{balsarahlle2d,balsarahllc2d,balsarahlle3d,BalsaraMultiDRS,MUSIC1,MUSIC2}; 
ii) the the second class uses divergence cleaning techniques, like either the Powell source term \cite{PowellMHD1} 
based  on the symmetric form of the MHD equations found by Godunov \cite{God1972MHD} or the hyperbolic generalized Lagrangian multiplier (GLM) approach of Munz et al. 
\cite{MunzCleaning} and Dedner et al. \cite{Dedneretal}. The method proposed in this paper falls into the first class of exactly divergence-free schemes. 
 
The rest of the paper is organized as follows: for the sake of simplicity and to facilitate the reader, we first present our new algorithm only for the ideal MHD equations 
in one space dimension, see Section \ref{sec.method}. Computational results for the one-dimensional case are shown in Section \ref{sec.results}. 
The extension of the method to the two-dimensional case, including viscous and resistive effects and a divergence-free evolution of the magnetic field is presented 
in Section \ref{sec.multid}. A set of classical benchmark problems for the ideal and viscous and resistive MHD equations is then solved in Section \ref{sec.results2d}, 
showing the performance of the method in the low Mach number limit as well as its robustness for shocked flows. 
Finally, in Section \ref{sec.concl} we give some concluding remarks and an outlook to future work.

\section{Numerical method for the ideal MHD equations in one space dimension} 
\label{sec.method}

\subsection{Governing PDE} 
\label{sec.mhd} 

The ideal MHD equations in one space dimension read as follows:   
\renewcommand{\arraystretch}{1.2} 
\begin{eqnarray}
\label{eqn.pde1d}    
		\frac{\partial}{\partial t} \left( \begin{array}{c} \rho \\ \rho u \\ \rho v \\ \rho w \\ \rho E \\ B_x \\ B_y \\ B_z \end{array}  \right) + 
		\frac{\partial}{\partial x} \left( 
		\begin{array}{c} 
		\rho u \\ 
		\rho u^2 + p + \frac{1}{8 \pi} \mathbf{B}^2 - \frac{1}{4 \pi} B_x^2   \\ 
		\rho u v           - \frac{1}{4 \pi} B_x B_y \\ 
		\rho u w           - \frac{1}{4 \pi} B_x B_z \\ 
		u \left( \rho E + p + \frac{1}{8 \pi} \mathbf{B}^2  \right) - \frac{1}{4 \pi} B_x (\mathbf{v}\cdot\mathbf{B})   \\ 
		0 \\ 
		u B_y - v B_x \\ 
		u B_z - w B_x 
		\end{array}  \right)                      & = & 0.  
\end{eqnarray} 
Here, time is denoted by $t \in \mathds{R}_0^+$, while $x \in \Omega = [x_L, x_R] \subset \mathds{R}$ is the spatial coordinate within the computational domain $\Omega$. 
As usual, the fluid density and the fluid pressure are denoted by $\rho$ and $p$, respectively;  $\mathbf{v}=(u,v,w)$ is the velocity field and the magnetic field
vector is $\mathbf{B}=(B_x,B_y,B_z)$; the total energy density is given by 
$\rho E = \rho e + \rho k + m = \rho e + \frac{1}{2} \rho \mathbf{v}^2 + \frac{1}{8\pi} \mathbf{B}^2$, where 
$\rho k = \frac{1}{2} \rho \mathbf{v}^2$ is the kinetic energy density of the fluid and $m=\frac{1}{8\pi} \mathbf{B}^2$ is the magnetic energy density; 
$e=e(p,\rho)$ is the specific internal energy per unit mass given by the so-called equation of state (EOS), which is in general a nonlinear function of the 
fluid pressure and density. However, for an ideal gas, $e$ is a linear function in $p$. In density-based Godunov-type finite volume schemes the EOS is typically required
in the form $p=p(e,\rho)$, which can be obtained by solving the expression $e=e(p,\rho)$ for the pressure. Another important quantity that we will use in this paper is 
the so-called specific enthalpy, which is defined as $h=e + p/\rho$ and which allows to rewrite the first part of the flux for the total energy density as follows: 
$u(\rho E + p + m)  = u (\rho k + 2 m) + h (\rho u)$. 
The eight eigenvalues of the MHD system \eqref{eqn.pde1d} are 
\begin{equation}
\lambda_{1,8} = u \mp c_f, \quad \lambda_{2,7} = u \mp c_a, \quad \lambda_{3,6} = u \mp c_s, \quad \lambda_4 = u, \quad \lambda_5 = 0, 
\label{eqn.eval.full} 
\end{equation} 
with 
\begin{equation}
 c_a = B_x / \sqrt{4 \pi \rho}, \quad
 c_s^2 = \halb \left( b+c - \sqrt{(b+c)^2-4 b_x c} \right), \quad 
 c_f^2 = \halb \left( b+c + \sqrt{(b+c)^2-4 b_x c } \right).  
\label{eqn.mhd.wavespeeds} 
\end{equation} 
Here, $c_a$ is the Alfv\'en wave speed, $c_s$ is the speed of the slow magnetosonic waves, $c_f$ is the one of the fast magnetosonic waves
and $c$ is the adiabatic sound speed that can be computed from the equation of state $p=p(e,\rho)$ as $c^2 = \partial p / \partial \rho + p / \rho^2 \partial p / \partial e$,
which reduces to the well-known expression $c^2 = \gamma p/ \rho$ for the ideal gas EOS. 
In the previous expressions we have also used the abbreviations $ b^2 = \mathbf{B}^2/ (4 \pi \rho)$ and $b_x = B_x / \sqrt{4 \pi \rho}$. 
Similar to the compressible Euler equations, the crucial terms that give rise to the fast and slow magnetosonic wave speed $c_f$ and $c_s$ are the pressure term 
$p_x$ in the momentum equation and the enthalpy term $( h \rho u)_x$ in the total energy equation. 
Therefore, these terms will have to be discretized \textit{implicitly} in our semi-implicit numerical method in order to avoid a CFL condition based on the magnetosonic 
wave speeds $c_f$ and $c_s$, while all remaining terms do not include the pressure and can therefore be discretized explicitly. 
In subsection \ref{sec.split}, we present a detailed discussion of the eigenvalues of an appropriately split MHD system, in order to properly motivate our choice for discretizing 
certain terms explicitly and others implicitly.  For a detailed analysis in the case of the compressible Euler and the shallow water equations, see \cite{CasulliCompressible,CasulliCattani}.  

\subsection{Ideal gas EOS} 
\label{sec.eos} 

Our numerical scheme is presented for a general nonlinear equation of state $e=e(p,\rho)$. However, in order to compare with previously published results in the literature, 
we will use the ideal gas EOS for all numerical test problems reported later. The ideal gas EOS in the sought form $e=e(p,\rho)$ can be obtained from the so-called 
thermal equation of state $p=p(\rho,T)$ and the so-called caloric equation of state $e=e(T,\rho)$ by eliminating the temperature. For the ideal gas, the thermal and caloric
equations of state take the well-known form 
\begin{equation} 
\label{eqn.thermcal.ideal} 
 \frac{p}{\rho} = R T, \qquad \textnormal{ and } \qquad e = c_v T,  
\end{equation} 
with the specific gas constant $R = c_p - c_v$, and the heat capacities $c_v$ and $c_p$ at constant volume and at constant pressure, respectively. From 
\eqref{eqn.thermcal.ideal} one easily obtains 
\begin{equation}
 e = e(p,\rho) = \frac{p}{(\gamma-1) \rho} ,  
\end{equation} 
which is linear in the pressure $p$ and where $\gamma = c_p / c_v$ denotes the so-called ratio of specific heats. For more general cubic equations 
of state, the reader is referred to the famous work by van der Waals \cite{vanderWaals} and more recent extensions, see 
\cite{Vidal,PengRobinson,RedlichKwong}. For completely general equations of state for real fluids, see \cite{Wagneretal,WagnerPruss}. 

\subsection{Split form of the MHD system} 
\label{sec.split} 

Following the seminal paper of Toro and V\'azquez \cite{ToroVazquez} on the Euler equations we now decide to \textit{split} the flux of the MHD 
system \eqref{eqn.pde1d} into a \textit{convective-type} flux and a \textit{pure pressure flux}, where the convective-type flux has to be understood in a more general sense 
in the MHD context due to the presence of the Alfv\'en waves. Note that the Toro \& V\'azquez (TV) splitting is substantially different from the flux vector 
splittings proposed in \cite{StegerWarming,ZhaBilgen,LiouSteffen}, since only in the TV splitting the resulting convective flux is totally free of any
pressure terms. Writing PDE \eqref{eqn.pde1d} formally as 
\begin{equation}
   \frac{\partial \Q}{\partial t} + \frac{\partial \f}{\partial x} = 0, 
	\label{eqn.pde.full} 
\end{equation} 
with $\Q = \left( \rho, \rho \mathbf{v}, \rho E, \mathbf{B} \right)$ the vector of conservative variables and the flux vector $\f$ given in \eqref{eqn.pde1d}, we write
the split system now as 
\begin{equation}
   \frac{\partial \Q}{\partial t} + \frac{\partial \f^c}{\partial x} + \frac{\partial \f^p}{\partial x} = 0, 
	\label{eqn.pde.split} 
\end{equation}
with the convective-type flux $\f^c$ and the pure pressure flux $\f^p$ given as follows: 
\begin{equation}
\f^c(\Q) = \left( 
		\begin{array}{c} 
		\rho u \\ 
		\rho u^2 + m - \frac{1}{4 \pi} B_x^2   \\ 
		\rho u v           - \frac{1}{4 \pi} B_x B_y \\ 
		\rho u w           - \frac{1}{4 \pi} B_x B_z \\ 
		u \left( \rho k + 2 m  \right) - \frac{1}{4 \pi} B_x (\mathbf{v}\cdot\mathbf{B})   \\ 
		0 \\ 
		u B_y - v B_x \\ 
		u B_z - w B_x 
		\end{array}  \right), \qquad 
\f^p(\Q) = \left( 		\begin{array}{c} 
		0   \\ 
		p   \\ 
		0   \\ 
		0   \\ 
		h \rho u  \\ 
		0   \\ 
		0   \\ 
		0     
		\end{array}  \right). 
\label{eqn.splitflux} 
\end{equation} 
Recall that we have $h=e(p,\rho)+p/\rho$ and $m= \frac{1}{8 \pi} \mathbf{B}^2$ as well as $\rho k = \frac{1}{2} \rho \mathbf{v}^2$ according to the definitions in subsection \ref{sec.mhd}. 
It is obvious that $\f^c$ does not contain any contribution of the pressure $p$, while $\f^p$ involves \textit{only} the density $\rho$, the pressure $p$ and the 
velocity component $u$ and does \textit{not} involve any contribution from the magnetic field. In this sense, our new splitting is the closest possible to the TV splitting, since our pressure system 
is exactly the same as the one obtained by Toro \& V\'azquez \cite{ToroVazquez}. Note in particular also that the split form \eqref{eqn.pde.split} with \eqref{eqn.splitflux} chosen in this paper is \textit{different} 
from all splittings of the MHD system proposed in Balsara et al. \cite{BalsaraMHDSplitting}. It is easy to check that the convective subsystem 
\begin{equation}
 \partial_t \Q + \partial_x \f^c = 0
\label{eqn.convective.sys} 
\end{equation} 
has the following eigenvalues 
\begin{equation}
 \lambda^c_{1,8} = u \mp \sqrt{\frac{\mathbf{B}^2}{4 \pi \rho}}, \qquad  
 \lambda^c_{2,7} = u \mp \frac{B_x}{\sqrt{ 4 \pi \rho }}, \qquad 
 \lambda^c_{3,4} = 0, \qquad 
 \lambda^c_{5,6} = u, 
	\label{eqn.eval.c} 
\end{equation}
while our pressure subsystem 
\begin{equation}
 \partial_t \Q + \partial_x \f^p = 0
\label{eqn.pressure.sys} 
\end{equation} 
is \textit{identical} to the Toro \& V\'azquez pressure system and therefore has the eigenvalues 
\begin{equation}
   \lambda^p_1 = \frac{1}{2} \left( u - \sqrt{u^2 + 4 c^2 } \right), \qquad \lambda^p_{2,3,4,5,6,7} = 0, \qquad 
	 \lambda^p_8 = \frac{1}{2} \left( u + \sqrt{u^2 + 4 c^2 } \right),   
	\label{eqn.eval.p} 
\end{equation} 
i.e. the pressure subsystem is \textit{always} subsonic, independent of $\mathbf{B}$. Recall that for the ideal gas EOS we have $c^2 = \gamma p/\rho$. 
Looking at the eigenvalues \eqref{eqn.eval.p} of the pressure subsystem it becomes obvious that for low Mach number flows, i.e. when the ratio 
$ M = u/c \ll 1$, or even more in the incompressible limit when $M \to 0$, the terms appearing in the pressure subsystem need to be discretized \textit{implicitly},  
while the eigenvalues of the convective subsystem \eqref{eqn.eval.c} suggest that an \textit{explicit} discretization of the convective subsystem is still possible 
unless the magnitude of the magnetic field $|\mathbf{B}|$ gets very large or the density gets very low. 
In the latter case, also the magnetic field needs to be discretized implicitly, but this is not the scope of the present paper. 

\subsection{Semi-implicit discretization on a staggered grid} 

The ideal MHD equations \eqref{eqn.pde1d} are discretized on a \textit{staggered} grid, which is typical for semi-implicit schemes 
applied to the incompressible Navier-Stokes and shallow water equations, see \cite{markerandcell,Casulli1990}. In the staggered grid the \textit{primary} control volumes  
are the intervals $\Omega_i=[x_{i-\halb},x_{i+\halb}]$ of length $\Delta x_i = x_{i+\halb}-x_{i-\halb}$ with barycenters located in $x_i = \halb ( x_{i-\halb} + x_{i+\halb} )$. 
The number of primary control volumes is denoted by $N_x$. The $N_x+1$ \textit{dual} control volumes are $\Omega_{i+\halb}=[x_i,x_{i+1}]$ with the associated mesh spacing 
$\Delta x_{i+\halb} = x_{i+1}-x_i = \halb (\Delta x_{i} + \Delta x_{i+1})$. 
The entire nonlinear convective subsystem \eqref{eqn.convective.sys} will be discretized on the primary control volumes, while the pressure subsystem \eqref{eqn.pressure.sys} 
is discretized as usual on the combination of the two staggered grids, defining the discrete pressure $p_i^n$ in the centers of the primary cells $\Omega_i$, 
while the discrete velocity $u_{i+\halb}^n$ in the pressure system is located at the cell boundaries. In order to combine the discretization of the 
convective subsystem on the main grid with the discretization of the pressure subsystem on the staggered mesh, we will need to \textit{average} quantities from the  
main grid to the dual grid, and vice versa. This is simply obtained by the following conservative averaging operators 
\begin{equation}
   \Q_i^n = \halb \left( \Q_{i-\halb}^n + \Q_{i+\halb}^n \right), \qquad 
	\Q_{i+\halb}^n = \halb \frac{1}{\Delta x_{i+\halb}} \left( \Delta x_i \Q_{i}^n + \Delta x_{i+1} \Q_{i+1}^n \right). 
	\label{eqn.average} 
\end{equation} 

\subsubsection{Convective subsystem.} 
The convective terms collected in $\f^c$ are now discretized on the main grid using a standard \textit{explicit} first or second order accurate finite volume scheme of the form
\begin{equation}
  \Q_i^{*} = \Q_i^n - \frac{\Delta t}{\Delta x_i}\left( \f^c_{i+\halb} - \f^c_{i-\halb} \right),
	\label{eqn.exp.fv} 
\end{equation} 
which yields the intermediate state vector $\Q_i^*$ that does not yet contain the contribution of the pressure terms. Throughout this paper we employ the simple Rusanov-type flux 
\begin{equation}
   \f^c_{i+\halb} = \halb \left( \f^c(\Q_{i+\halb}^+) + \f^c(\Q_{i-\halb}^+) \right) - \halb s_{\max} \left( \Q_{i+\halb}^+ - \Q_{i+\halb}^- \right),  
	\label{eqn.rusanov}
\end{equation} 
where $\Q_{i+\halb}^-$ and $\Q_{i+\halb}^+$ denote the left and right boundary extrapolated states at the cell interface $x_{i+\halb}$ and  
$s_{\max} = \max \left( |\lambda^c_l(\Q_{i+\halb}^-) |, |\lambda^c_l(\Q_{i+\halb}^+) | \right)$ is the maximum signal speed of the convective subsystem at the interface. 
For a first order scheme one simply has $\Q_{i+\halb}^- = \Q_i^n$ and $\Q_{i+\halb}^+ = \Q_{i+1}^n$, while a second order MUSCL-Hancock-type TVD scheme is obtained by setting 
$\Q_{i+\halb}^- = \mathbf{w}_i(x_{i+\halb},t^{n+\halb})$ and $\Q_{i+\halb}^+ = \mathbf{w}_{i+1}(x_{i+\halb},t^{n+\halb})$, where $\mathbf{w}_i(x,t)$ is a space-time 
polynomial reconstruction of the state vector in each cell $\Omega_i$ that reads 
\begin{equation}
  \mathbf{w}_i(x,t) = \Q_i^n + \frac{\Delta \Q_i^n}{\Delta x_i}\left( x - x_i \right) + \partial_t \Q_i^n \left( t - t^n \right). 
	\label{eqn.rec} 
\end{equation} 
The space-time expansion coefficients in \eqref{eqn.rec} are given by 
\begin{equation}
   \frac{\Delta \Q_i^n}{\Delta x_i} = \textnormal{minmod}\left( \frac{\Q_{i+1}^n - \Q_i^n}{\Delta x_{i+\halb}}, \frac{\Q_i^n - \Q_{i-1}^n}{\Delta x_{i-\halb}} \right), \,  
		\partial_t \Q_i^n = \frac{\f^c \left( \Q_i^n - \halb \Delta \Q_i^n \right) - \f^c \left( \Q_i^n + \halb \Delta \Q_i^n \right)}{\Delta x_i},
\end{equation} 
with the usual minmod slope limiter function, see \cite{toro-book}. Since the mass conservation equation and the PDE for the transverse momentum in $y$ and $z$ direction do not 
contain the pressure, we can immediately set $\rho_i^{n+1}=\rho_i^*$, $(\rho v)_i^{n+1}=(\rho v)_i^{*}$ and $(\rho w)_i^{n+1}=(\rho w)_i^{*}$. 
In one space dimension, no divergence-free treatment of the magnetic field is necessary and therefore we also have $\mathbf{B}_i^{n+1} = \mathbf{B}_i^*$. 
This completes the description of the explicit part of the scheme. 

\subsubsection{Pressure subsystem.} 
The pressure subsystem involves only the $x$-momentum equation and the total energy equation. The \textit{semi-implicit} discretization of the $x$-momentum equation reads 
\begin{equation}
\label{eqn.rhou.disc} 
 (\rho u)_{i+\halb}^{n+1} = (\rho u)_{i+\halb}^{*} - \frac{\Delta t}{\Delta x_{i+\halb}} \left( p_{i+1}^{n+1} - p_i^{n+1} \right),  
\end{equation} 
where the pressure is now taken \textit{implicitly}, while the \textit{explicit} operator for the discretization of the nonlinear convective terms i.e. 
for the computation of $(\rho u)_{i+\halb}^{*}$ has been detailed previously. Note that $(\rho u)_{i+\halb}^{*}$ is located on the dual mesh and therefore
has to be averaged from the main grid to the dual mesh via \eqref{eqn.average}.  
According to \cite{DumbserCasulli2016} a \textit{preliminary} discretization of the total energy equation is now chosen as follows: 
\begin{eqnarray}
\label{eqn.rhoE.prelim} 
 \Delta x_{i} \left( \rho_i^{n+1} e\left( p_i^{n+1}, \rho_i^{n+1} \right) +  \frac{1}{2} \left( (\tilde{\rho k})^{n+1}_{i-\halb} + (\tilde{\rho k})^{n+1}_{i+\halb} \right) 
+ m_i^{n+1} \right) = &&   \nonumber \\ 
\Delta x_{i} (\rho E)_i^* - \Delta t \left( \tilde{h}_{i+\halb}^{n+1} (\rho u)_{i+\halb}^{n+1} - \tilde{h}_{i-\halb}^{n+1} (\rho u)_{i-\halb}^{n+1} \right)  .  && 
\end{eqnarray} 
The tilde symbols indicate that a further discretization step is necessary that will be explained later. Inserting the discrete momentum equation \eqref{eqn.rhou.disc} into 
the discrete energy equation \eqref{eqn.rhoE.prelim} and using $m_i^{n+1}=m_i^*$ yields the following preliminary system for the unknown pressure 
\begin{eqnarray}
\label{eqn.p.prelim} 
 \Delta x_{i} \rho_i^{n+1} e\left( p_i^{n+1}, \rho_i^{n+1} \right) - {\Delta t^2} \left( \frac{\tilde{h}_{i+\halb}^{n+1}}{\Delta x_{i+\halb}} \left( p_{i+1}^{n+1}-p_{i}^{n+1} \right) - 
\frac{\tilde{h}_{i-\halb}^{n+1}}{\Delta x_{i-\halb}} \left( p_{i}^{n+1}-p_{i-1}^{n+1} \right) \right) & = & \nonumber \\  
\Delta x_{i} \left( (\rho E)_i^* - m_i^* - \frac{1}{2} \left( (\tilde{\rho k})^{n+1}_{i-\halb} + (\tilde{\rho k})^{n+1}_{i+\halb} \right) \right) 
              - {\Delta t} \left( \tilde{h}_{i+\halb}^{n+1} (\rho u)_{i+\halb}^{*} - \tilde{h}_{i-\halb}^{n+1} (\rho u)_{i-\halb}^{*} \right),  &&  
\end{eqnarray} 
which has exactly the same structure as the one obtained in \cite{DumbserCasulli2016} for the compressible Euler equations. Therefore, following the same reasoning as explained
in \cite{DumbserCasulli2016}, the quantities marked with a tilde symbol cannot be discretized directly at the new time $t^{n+1}$, since in this case the resulting pressure system 
would become strongly nonlinear and difficult to control. To circumvent the problem, we employ a simple \textit{Picard iteration}, as suggested in \cite{CasulliZanolli2010}. The 
Picard iteration index will be denoted by $r$ in the following. This yields the following iterative scheme, which requires only the solution of the following \textit{mildly nonlinear} 
system for the pressure $p_i^{n+1,r+1}$ at each Picard iteration: 
\begin{equation}
\label{eqn.p} 
 \Delta x_{i} \rho e\left( p_i^{n+1,r+1} \right) - {\Delta t^2} \! \left( \! \frac{{h}_{i+\halb}^{n+1,r}}{\Delta x_{i+\halb}} \left( p_{i+1}^{n+1,r+1} \! -p_{i}^{n+1,r+1} \! \right) - 
\frac{{h}_{i-\halb}^{n+1,r}}{\Delta x_{i-\halb}} \left( p_{i}^{n+1,r+1} \! -p_{i-1}^{n+1,r+1} \right) \! \right) \! = b_i^r, 
\end{equation} 
with the abbreviation $\rho e\left( p_i^{n+1,r+1} \right) =  \rho_i^{n+1} e\left( p_i^{n+1,r+1}, \rho_i^{n+1} \right)$ and the known right hand side 
\begin{equation} 
b_i^r = \Delta x_{i} \! \left( (\rho E)_i^* - m_i^* - \frac{1}{2} \left( ({\rho k})^{n+1,r}_{i-\halb} + ({\rho k})^{n+1,r}_{i+\halb} \right) \right) \! 
              - \Delta t \left( {h}_{i+\halb}^{n+1,r} (\rho u)_{i+\halb}^{*} - {h}_{i-\halb}^{n+1,r} (\rho u)_{i-\halb}^{*} \right). 
\end{equation} 
Note that the density $\rho_i^{n+1}=\rho_i^{*}$ and the magnetic energy $m_i^{n+1}=m_i^*$ are already known from the explicit discretization \eqref{eqn.exp.fv}, hence in \eqref{eqn.p} the new pressure 
is the only unknown. Using a more compact notation, the above system \eqref{eqn.p} can be written as follows:  
\begin{equation}
\label{eqn.nonlinear} 
 \boldsymbol{\rho}\mathbf{e}(\mathbf{p}^{n+1,r+1}) + \mathbf{T}^r \, \mathbf{p}^{n+1,r+1} = \mathbf{b}^r, 
\end{equation} 
with the vector of the unknowns $\mathbf{p}^{n+1,r+1}=(p_1^{n+1,r+1}, ..., p_i^{n+1,r+1}, ..., p_{N_x}^{n+1,r+1})$. 
The vector $\mathbf{b}^r$ contains the known right hand side of \eqref{eqn.p}. 
Matrix $\mathbf{T}^r$ is symmetric and at least positive semi-definite and takes into account the linear part of the system, while the nonlinearity 
is contained in the vector function $\boldsymbol{\rho}\mathbf{e}(\mathbf{p}^{n+1,r+1}) = \left( \Delta x_1 \rho_1^{n+1} e(p_1^{n+1,r+1},\rho_1^{n+1}), ..., \Delta x_i \rho_i^{n+1} e(p_i^{n+1,r+1},\rho_i^{n+1}), ..., 
\Delta x_{N_x} \rho_{N_x}^{n+1} e(p_{N_x}^{n+1,r+1},\rho_{N_x}^{n+1})\right)$, 
which means a componentwise evaluation of the internal energy density in terms of pressure and density. We stress again that the density $\rho_i^{n+1}$ at the new time level is already 
known from \eqref{eqn.exp.fv}, i.e. for the solution of the mildly nonlinear system, the equation of state can be considered as a function of pressure alone, with a given density. 

The time step, the mesh spacings and the enthalpy $h$ are non-negative quantities and we suppose that the specific internal energy $e\left( p, \rho \right)$ is a non-negative, non-decreasing 
function whose derivative w.r.t. the pressure is a function of bounded variation. 
Thanks to the semi-implicit discretization of the pressure subsystem on the staggered mesh, the matrix $\mathbf{T}^r$ in system \eqref{eqn.p} is 
\textit{symmetric} and at least \textit{positive semi-definite}, which is quite a remarkable property, considering the complex structure of the MHD system \eqref{eqn.pde1d}. 
It is therefore possible to employ the same (nested) Newton-type techniques for the solution of \eqref{eqn.nonlinear} as those proposed and analyzed 
by Casulli et al. in \cite{CasulliZanolli2010,CasulliZanolli2012,BrugnanoCasulli,BrugnanoCasulli2}. For all implementation details and 
a rigorous convergence proof of the (nested) Newton method, the reader is referred to the above references. The iterative Newton-type techniques of Casulli et al. have already been used 
with great success as building block of semi-implicit finite volume schemes in different application contexts, see \cite{Casulli2009,CasulliStelling2011,CasulliDumbserToro2012,CasulliVOF,BoscheriDumbser,FambriDumbserCasulli,DumbserIbenIoriatti}. Due to the properties of $\mathbf{T}^r$, the linear sub-problems within the Newton-type algorithm can be solved at the aid of a 
matrix-free conjugate gradient method, or with the Thomas algorithm for tri-diagonal systems in the one-dimensional case. 
Note that for the ideal gas EOS the resulting system \eqref{eqn.nonlinear} becomes \textit{linear} in the pressure, hence one single Newton iteration is sufficient to 
solve \eqref{eqn.nonlinear}. 
From the new pressure $p_i^{n+1,r+1}$ the momentum density at the next Picard iteration can be obtained as 
\begin{equation}
 \label{eqn.rhou.picard} 
(\rho u)_{i+\halb}^{n+1,r+1} = (\rho u)_{i+\halb}^{*} - \frac{\Delta t}{\Delta x_{i+\halb}} \left( p_{i+1}^{n+1,r+1} - p_i^{n+1,r+1} \right).    						
\end{equation} 
The new pressure and momentum are both needed to update the enthalpies at the element interfaces as well as the kinetic energy contribution to the total energy at the new time level. 
As already observed in \cite{CasulliZanolli2010,DumbserCasulli2016} it is sufficient to carry out only very few Picard iterations to obtain a satisfactory solution. In all test problems presented
in this paper, we stop the Picard process after $r_{\max}=2$ iterations. 
At the end of the last Picard iteration, we set $p_i^{n+1}=p_i^{n+1,r+1}$, $(\rho u)_{i+\halb}^{n+1} = (\rho u)_{i+\halb}^{n+1,r+1}$, 
$h_{i+\halb}^{n+1} = h_{i+\halb}^{n+1,r+1}$ and update the total energy  density using the conservative formula  
\begin{equation}
\label{eqn.rhoE.disc} 
 (\rho E)_i^{n+1} = (\rho E)_i^* - \frac{\Delta t}{\Delta x_{i}} \left( {h}_{i+\halb}^{n+1} (\rho u)_{i+\halb}^{n+1} - {h}_{i-\halb}^{n+1} (\rho u)_{i-\halb}^{n+1} \right). 
\end{equation} 
Finally, in order to proceed with the next time step, we still need to average the momentum back from the staggered mesh to the main grid by using the averaging operator 
\eqref{eqn.average} from the dual mesh to the main grid. 

From \eqref{eqn.exp.fv}, \eqref{eqn.rhou.disc} and \eqref{eqn.rhoE.disc} it is obvious that the scheme is written in a conservative flux form for all conservation equations and the averaging operators
between main and dual grid are also conservative, hence the proposed method is locally and globally conservative for mass, momentum and total energy. Its stability is only restricted by a 
\textit{mild CFL condition} based on the eigenvalues of the convective subsystem $\lambda^c_l$, and is \textit{not} based on the speed of the magnetosonic waves $c_s$ and $c_f$. 
This makes the method particularly well suited for the discretization of \textit{low Mach number} flows. However, being a locally and globally conservative scheme, it is also able to handle 
flows with very strong shocks properly, as shown via several numerical test problems in the next section. 


\section{Numerical results in 1D}
\label{sec.results} 

\begin{figure}[!t]
\begin{center}
\includegraphics[width=0.6\textwidth]{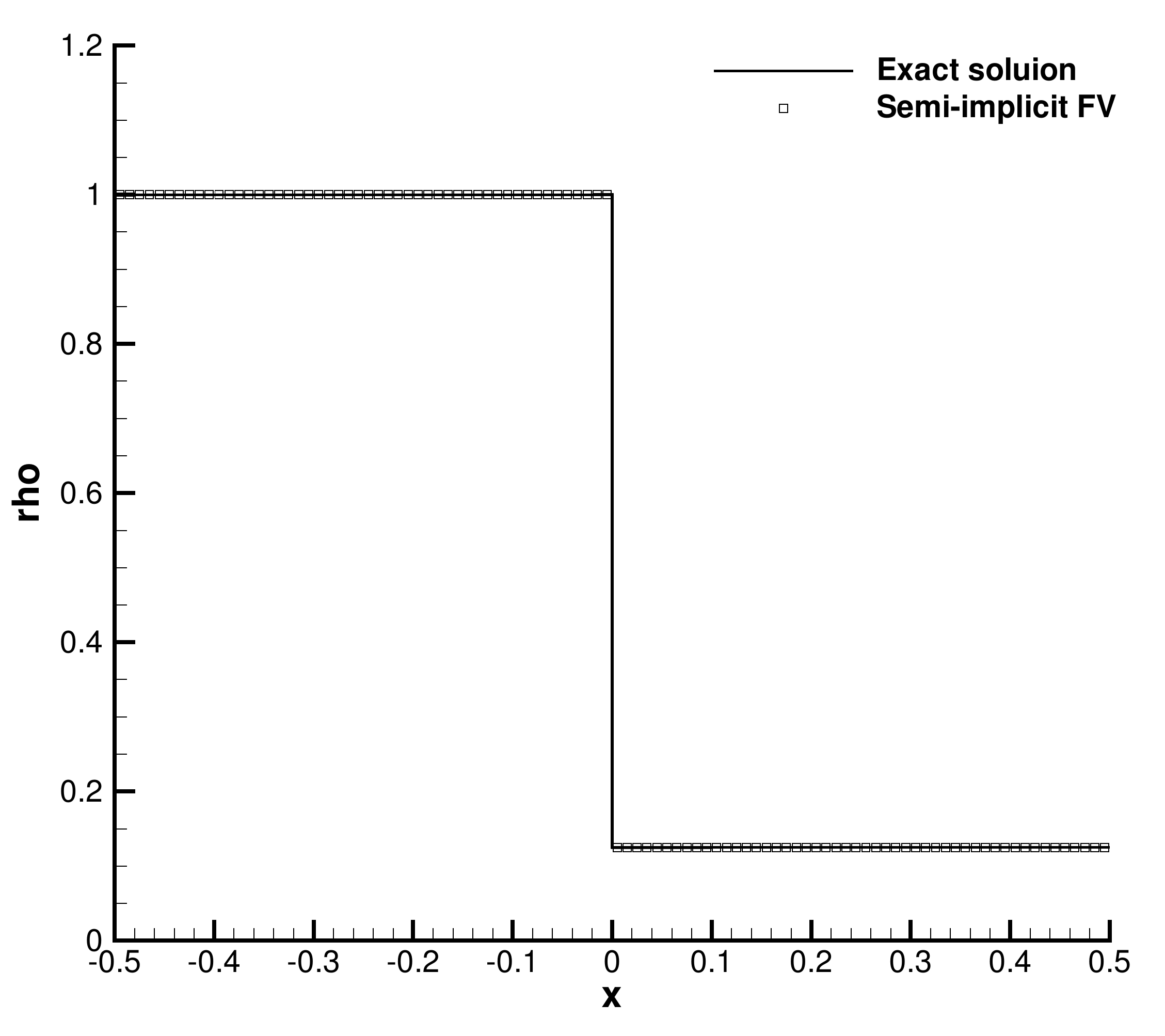}  
\caption{Exact and numerical solution for Riemann problem RP0 (isolated steady contact wave) solving the ideal MHD equations with the new SIFV scheme. 
The density is shown at a final time of $t=10$, confirming that our scheme is able to preserve steady contact waves exactly.} 
\label{fig.rp0}
\end{center}
\end{figure}

\begin{figure}[!t]
\begin{center}
\begin{tabular}{cc} 
\includegraphics[width=0.45\textwidth]{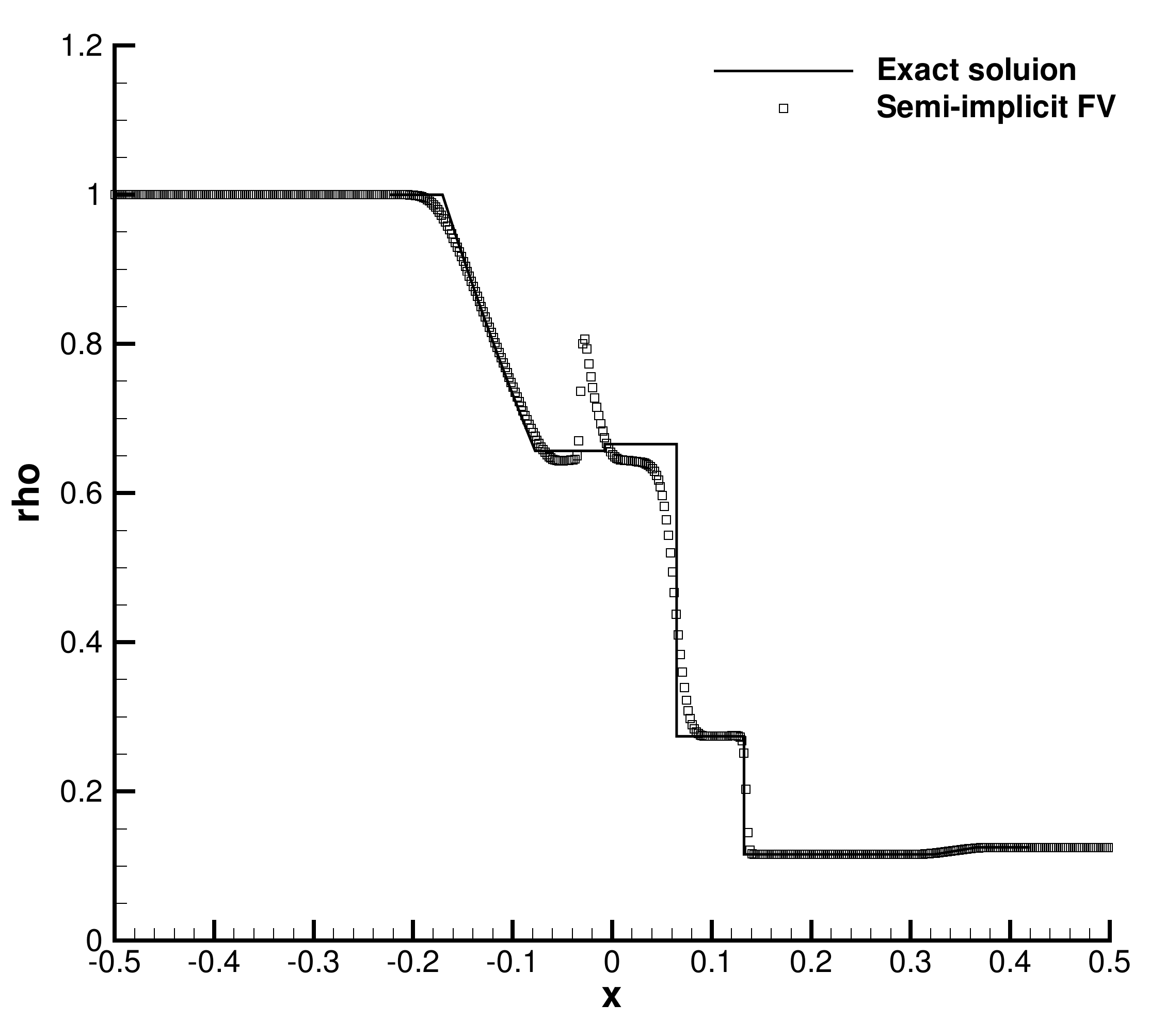} & 
\includegraphics[width=0.45\textwidth]{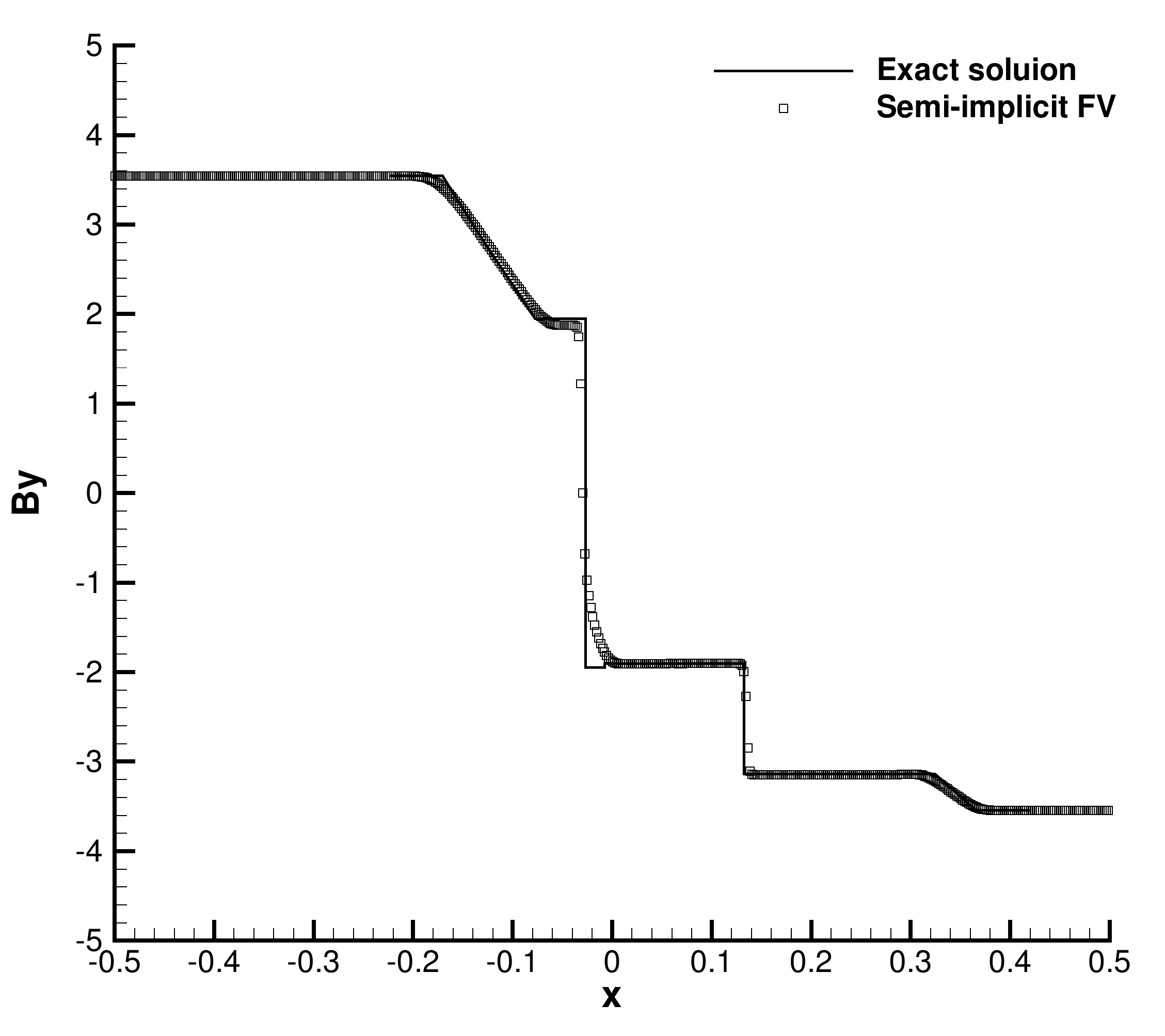}    
\end{tabular} 
\caption{Exact and numerical solution for Riemann problem RP1 solving the ideal MHD equations with the new SIFV scheme. Density (left) and magnetic field component $B_y$ (right) at time $t=0.1$. } 
\label{fig.rp1}
\end{center}
\end{figure}

\begin{figure}[!ht]
\begin{center}
\begin{tabular}{cc} 
\includegraphics[width=0.45\textwidth]{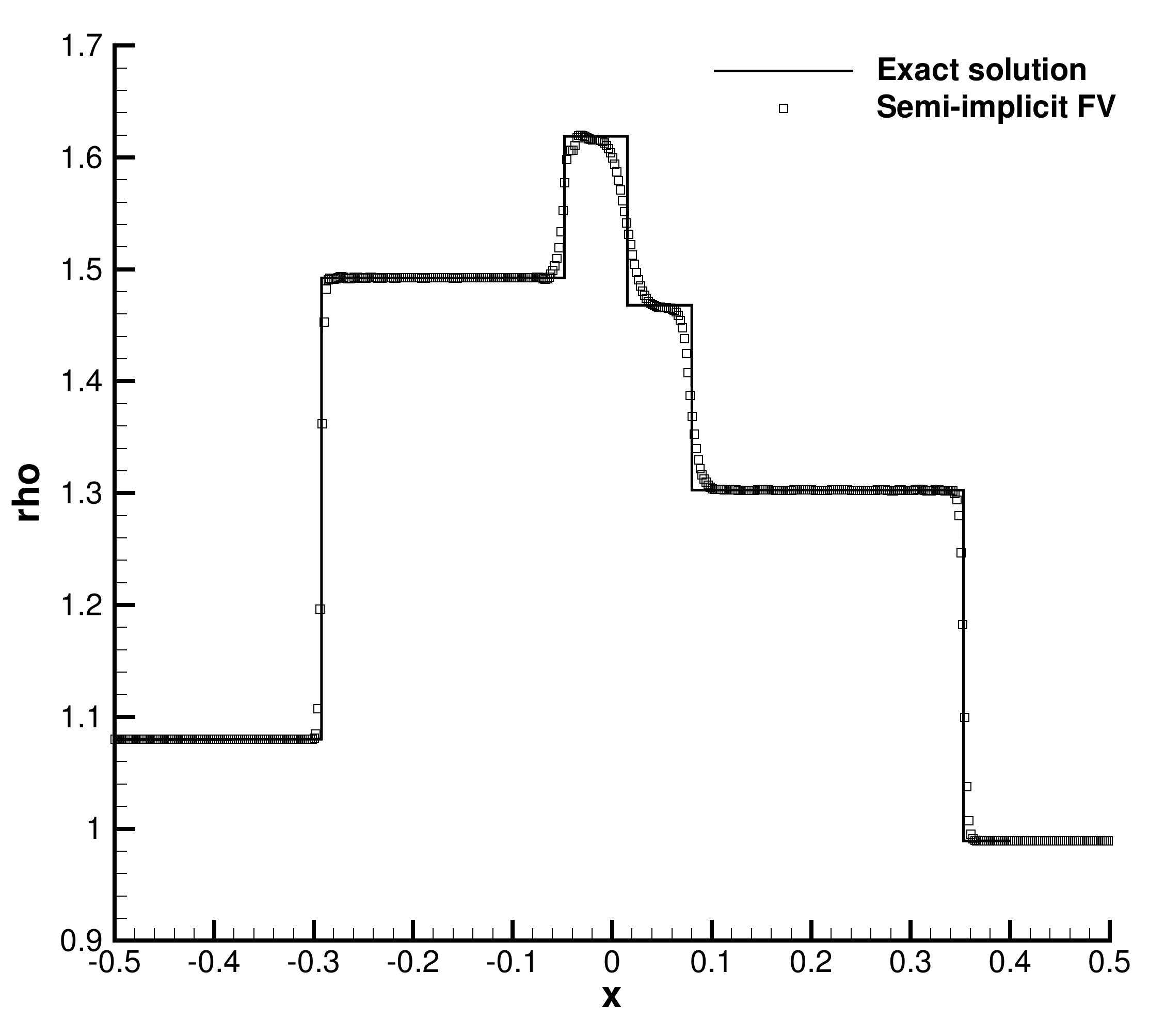} & 
\includegraphics[width=0.45\textwidth]{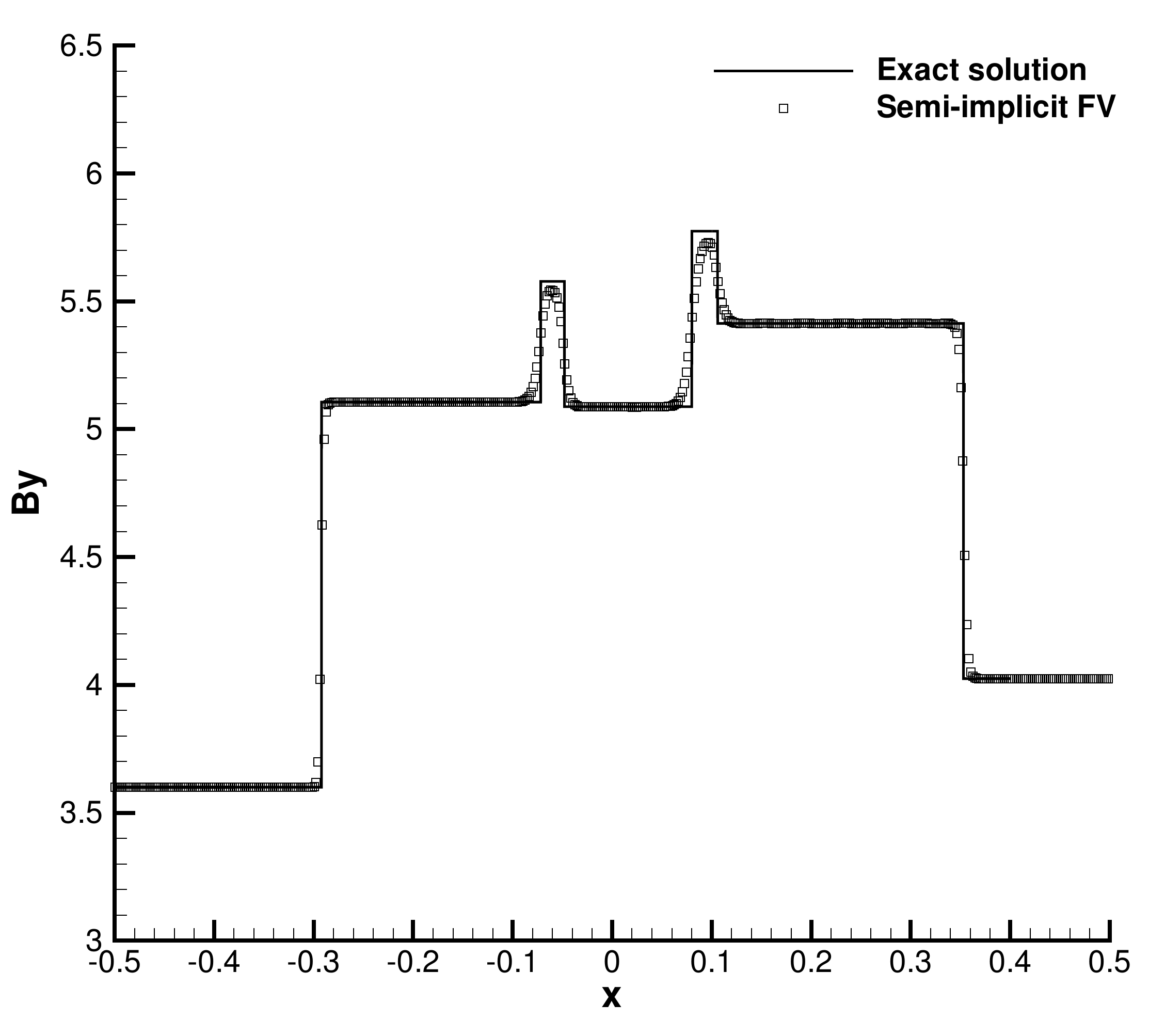}    
\end{tabular} 
\caption{Exact and numerical solution for Riemann problem RP2 solving the ideal MHD equations with the new SIFV scheme. Density (left) and magnetic field component $B_y$ (right) at time $t=0.2$. } 
\label{fig.rp2}
\end{center}
\end{figure}

\begin{figure}[!ht]
\begin{center}
\begin{tabular}{cc} 
\includegraphics[width=0.45\textwidth]{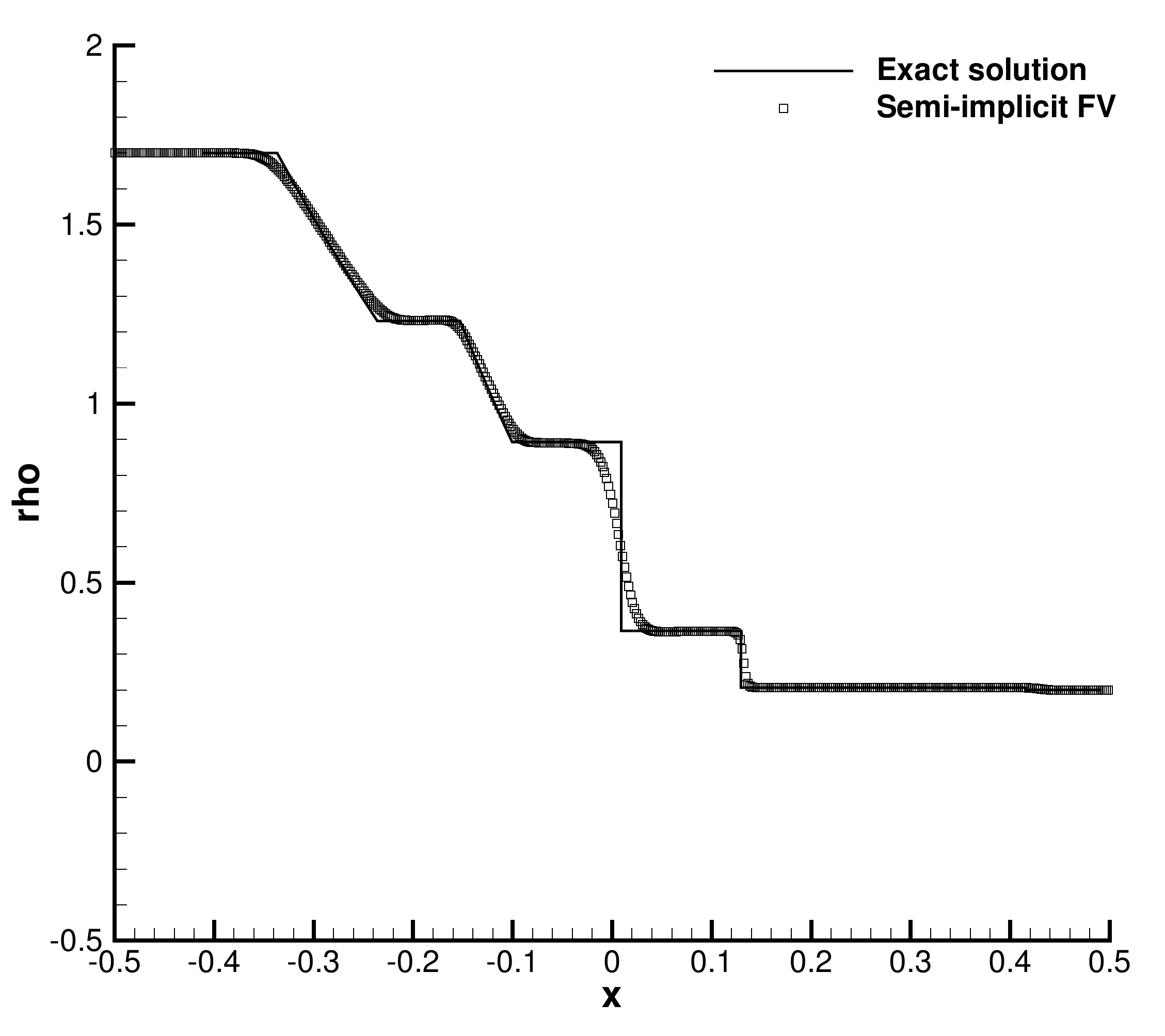} & 
\includegraphics[width=0.45\textwidth]{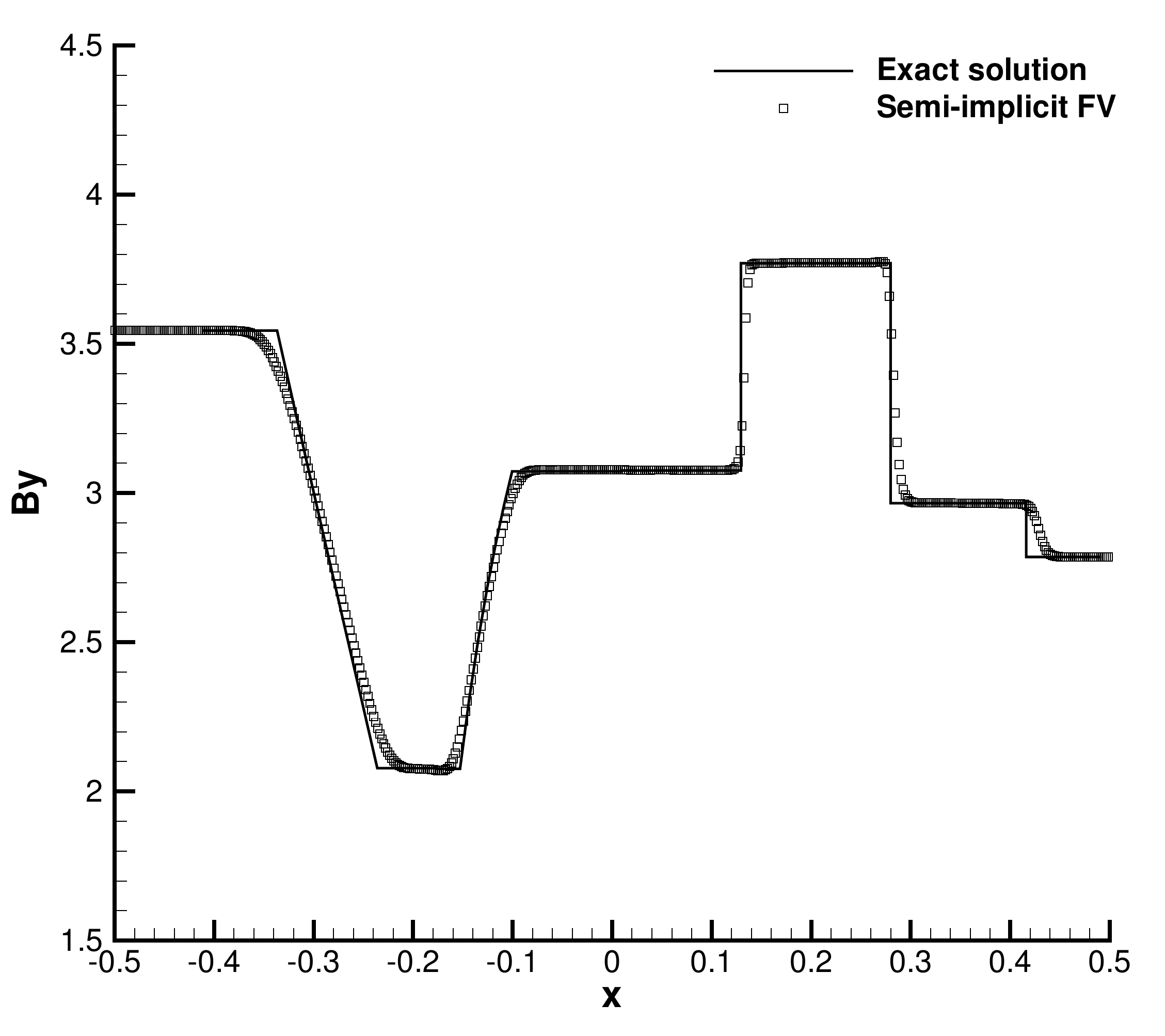}    
\end{tabular} 
\caption{Exact and numerical solution for Riemann problem RP3 solving the ideal MHD equations with the new SIFV scheme. Density (left) and magnetic field component $B_y$ (right) at time $t=0.15$. } 
\label{fig.rp3}
\end{center}
\end{figure}

\begin{figure}[!ht]
\begin{center}
\begin{tabular}{cc} 
\includegraphics[width=0.45\textwidth]{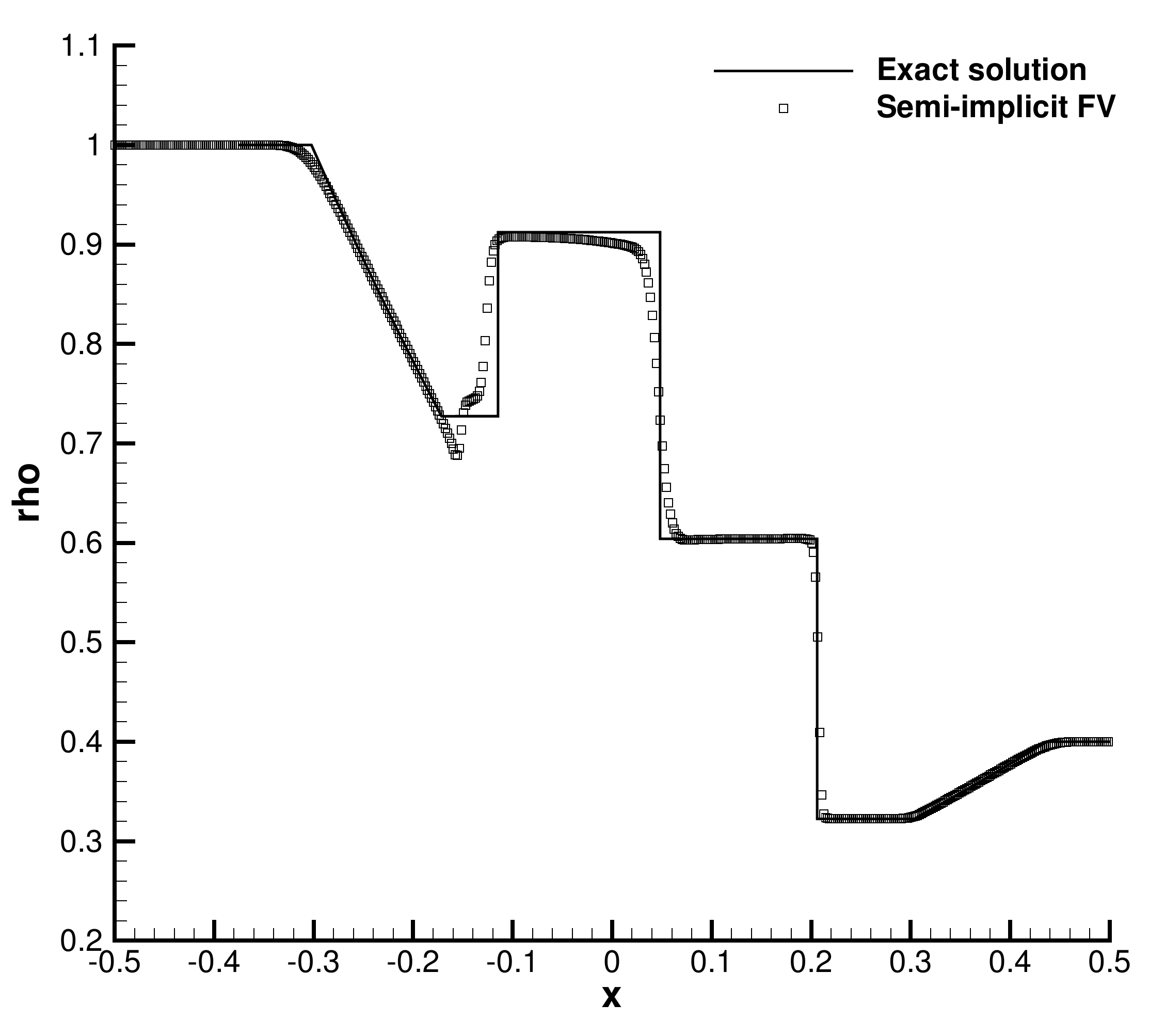} & 
\includegraphics[width=0.45\textwidth]{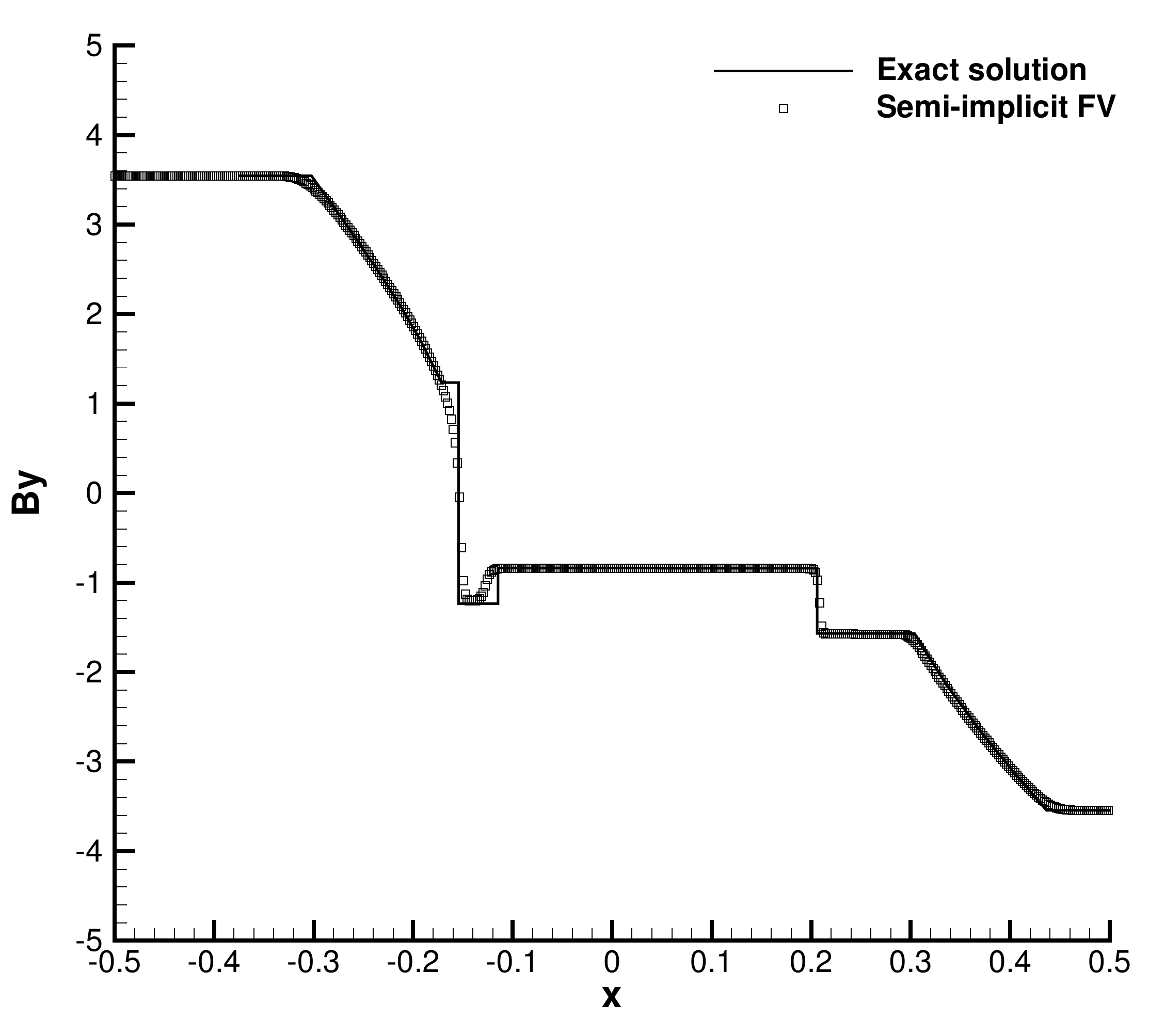}    
\end{tabular} 
\caption{Exact and numerical solution for Riemann problem RP4 solving the ideal MHD equations with the new SIFV scheme. Density (left) and magnetic field component $B_y$ (right) at time $t=0.16$. } 
\label{fig.rp4}
\end{center}
\end{figure}

In this section we apply our new semi-implicit finite volume scheme to a set of Riemann problems of the ideal MHD equations, some of which have been introduced 
and analyzed in \cite{BrioWu,RyuJones,DaiWoodward,OsherUniversal}. The eigenstructure of the ideal MHD equations has been discussed in \cite{roebalsara}, 
while the exact Riemann solver used for the comparisons presented in this paper has kindly been provided by S.A.E.G. Falle \cite{fallemhd,falle2}. 
For an alternative exact Riemann solver of the MHD equations see the work of Torrilhon \cite{torrilhon}. In all our tests we use a computational domain 
$\Omega=[-0.5,+0.5]$ that is discretized at the aid of 1000 pressure control volumes (apart from RP0, for which only 100 points have been used), 
which is only slightly more than the typical resolution of 800 elements chosen for the explicit finite volume schemes used in \cite{Balsara1998,torrilhon,torrilhonbalsara}. 
In RP1-RP4 the Courant number is set to $\CFL=0.9$, based on the maximum eigenvalues of the convective subsystem and a second order MUSCL-type TVD scheme is used for the 
discretization of the explicit terms. For RP0, we use a constant time step size of $\Delta t = 0.1$. 
The initial condition for all Riemann problems consists in a constant left and right state that are separated by a discontinuity located at $x_d$. The initial data 
as well as the value for $x_d$ are reported in Table \ref{tab.ic.mhd}. The ratio of specific heats is $\gamma=\frac{5}{3}$ for all cases. 
The comparison between the numerical solution obtained with the new SIFV scheme and the exact solution is presented in Figs. \ref{fig.rp0}-\ref{fig.rp4}. 
The first Riemann problem (RP0) is just a sanity check in order to verify that our new SIFV method is able to resolve isolated steady contact waves without magnetic field
exactly. This property follows trivially from the chosen discretization and is also confirmed in our numerical experiments, see Fig. \ref{fig.rp0}. 

Riemann problem (RP1) is the one of Brio \& Wu \cite{BrioWu}, for which it is well known that all standard finite volume schemes 
produce a compound wave instead of the wave pattern suggested by the exact Riemann solver. Only the random choice method of Glimm \cite{Glimm:1965a} was able to reproduce the 
correct solution in this case, as discussed in \cite{fallemhd}. Therefore, despite the disagreement with the exact solution in the density profile, 
our numerical results are in line with others published in the literature. Furthermore, the numerical results obtained for the magnetic field component $B_y$ agree well 
with the exact solution. 
The second Riemann problem (RP2) goes back to Ryu \& Jones \cite{RyuJones} and presents a wave pattern composed of discontinuities in all seven waves of the MHD system. 
The agreement between our numerical solution and the exact solution is very good in this case. Also problem RP3 contains seven waves, but compared to RP2 the two left waves 
are rarefactions and not shocks. Also in this case the semi-implicit finite volume scheme is able to capture the wave pattern properly, apart from the weak right-moving shock. 
In the last Riemann problem (RP4), our scheme has some difficulties in capturing the second wave from the left at the given grid resolution, but this behaviour is similar to 
what was also observed in \cite{OsherUniversal} and \cite{HPRmodelMHD}. The profile of the magnetic field component $B_y$ is well reproduced in this case. 

Overall we can conclude that the numerical results obtained with our new algorithm are in line with those previously published in the literature. However, at this point it is 
important to stress that our semi-implicit finite volume scheme is a so-called \textit{pressure-based} solver, which is particularly tailored to work in the low Mach number 
regime or even in the incompressible limit of the equations, while all standard explicit finite volume schemes that are typically used for the solution of the MHD equations 
are so-called \textit{density-based} methods, which are unable to deal with the incompressible limit of the equations. It is therefore quite remarkable to observe that the 
new pressure-based semi-implicit method performs almost as well as standard Godunov-type schemes in this set of Riemann problems. Encouraged by these results, in the next
section we now present the extension to the viscous and resistive case in two space dimensions, where particular care needs to be taken in order to obtain an exactly 
divergence-free formulation of the scheme. 

\begin{table}[!t]
 \caption{Initial states left and right for the density $\rho$, velocity vector $\mathbf{v} = (u,v,w)$, the pressure $p$  
 and the magnetic field vector $\mathbf{B} = (B_x,B_y,B_z)$ for the Riemann problems of the ideal classical MHD equations. 
 In all cases $\gamma=5/3$. The initial position of the discontinuity is $x_d=0$ for RP0, RP1 and RP4, while it is 
$x_d=-0.1$ for RP2 and RP3. } 
\begin{center} 
 \begin{tabular}{rcccccccc}
 \hline
 Case & $\rho$ & $u$ & $v$ & $w$ & $p$ & $B_x$ & $B_y$ & $B_z$        \\ 
 \hline   
 RP0 L: &  1.0    &  0.0     & 0.0    & 0.0      &  1.0     & 0.0 & 0.0  & 0.0       \\
     R: &  0.125  &  0.0     & 0.0    & 0.0      &  1.0     & 0.0 & 0.0  & 0.0       \\
 RP1 L: &  1.0    &  0.0     & 0.0    & 0.0      &  1.0     & $\frac{3}{4} \sqrt{4 \pi}$ &  $\sqrt{4 \pi}$  & 0.0       \\
     R: &  0.125  &  0.0     & 0.0    & 0.0      &  0.1     & $\frac{3}{4} \sqrt{4 \pi}$ & $-\sqrt{4 \pi}$  & 0.0       \\
 RP2 L: &  1.08   &  1.2     & 0.01   & 0.5      &  0.95    & 2.0 &  3.6     & 2.0            \\
     R: &  0.9891 &  -0.0131 & 0.0269 & 0.010037 &  0.97159 & 2.0 &  4.0244  & 2.0026         \\
 RP3 L: &  1.7    &  0.0     & 0.0    & 0.0      &  1.7     & 3.899398 &  3.544908  & 0.0              \\
     R: &  0.2    &  0.0     & 0.0    & -1.496891  &  0.2   & 3.899398 &  2.785898  & 2.192064         \\
 RP4 L: &  1.0    &  0.0     & 0.0    & 0.0      &  1.0     & $1.3 \sqrt{4 \pi}$ &  $\sqrt{4 \pi}$   & 0.0            \\
     R: &  0.4    &  0.0     & 0.0    & 0.0      &  0.4     & $1.3 \sqrt{4 \pi}$ &  $-\sqrt{4 \pi}$  & 0.0             \\
 \hline
 \end{tabular}
\end{center} 
 \label{tab.ic.mhd}
\end{table}

\section{Extension to viscous flows in multiple space dimensions}
\label{sec.multid} 

\subsection{Governing equations}
\label{sec.pde.multid} 

In multiple space dimensions, the viscous and resistive MHD equations read 
\begin{equation}
\label{eqn.pde2d}    
		\frac{\partial}{\partial t} \left( \begin{array}{c} \rho \\ \rho \mathbf{v} \\ \rho E \\ \mathbf{B} \end{array}  \right) 
		+ \nabla \cdot \left( \begin{array}{c} \rho \mathbf{v} \\ \rho \mathbf{v} \otimes \mathbf{v} +  \left( p + \frac{\mathbf{B}^2}{8\pi} \right) \mathbf{I} - \frac{1}{4\pi} \mathbf{B} \otimes \mathbf{B} \\  \mathbf{v}^T \left( \rho E + p + \frac{1}{8\pi} \mathbf{B}^2 \right)  - \frac{1}{4\pi} \mathbf{v}^T \mathbf{B} \otimes \mathbf{B} \\ \mathbf{B} \otimes \mathbf{v} - \mathbf{v} \otimes \mathbf{B} \end{array}  \right)   = 
		\nabla \cdot \mathbf{F}_v,  
\end{equation} 
with the identity matrix $\mathbf{I}$ and the viscous flux tensor $\mathbf{F}^v = (\mathbf{f}^v, \mathbf{g}^v) $ defined as 
\begin{equation}
\mathbf{F}^v(\mathbf{V}, \nabla \mathbf{V}) = \left( \begin{array}{c} 
0 \\ 
 \mu \left( \nabla \mathbf{v} + \nabla \mathbf{v}^T - \frac{2}{3} \left( \nabla \cdot \mathbf{v} \right) \mathbf{I} \right) \\
\mu \mathbf{v}^T  \left( \nabla \mathbf{v} + \nabla \mathbf{v}^T - \frac{2}{3} \left( \nabla \cdot \mathbf{v} \right) \mathbf{I} \right) + \lambda \nabla T + \frac{\eta}{4 \pi} \mathbf{B}^T \left( \nabla \mathbf{B} - \nabla \mathbf{B}^T \right) \\
\eta \left( \nabla \mathbf{B} - \nabla \mathbf{B}^T \right)    
\end{array} \right). 
\end{equation} 
Here $\mathbf{V}=(\rho,\mathbf{v},T,\mathbf{B})$ is the vector of primitive variables, $T$ is the temperature given by the thermal equation of state $T=T(p,\rho)$, 
$\mu$ is the kinematic viscosity, $\lambda$ is the thermal conductivity and $\eta$ is the electric resistivity of the fluid. The Prandtl number is defined as 
$Pr = \mu \gamma c_v / \lambda$. 
In order to extend our new semi-implicit scheme to the ideal and to the viscous and resistive MHD equations in multiple space dimensions, special care must be taken 
concerning the $\nabla \cdot \mathbf{B} = 0$ constraint, i.e. the divergence of the magnetic field must remain zero for all times if it was initially zero. Several 
strategies have been developed in the literature in the past in order to satisfy the divergence constraint exactly or approximately, see e.g. the well-known 
divergence-free schemes for MHD of Balsara and Spicer \cite{BalsaraSpicer1999} and subsequent work by Balsara \cite{Balsara2004,balsarahlle2d,balsarahllc2d}, 
the discretization proposed by Powell \cite{PowellMHD1,PowellMHD2} based on the symmetric hyperbolic form of the MHD equations found by Godunov \cite{God1972MHD}, 
or the hyperbolic divergence-cleaning approach of Munz et al. \cite{MunzCleaning} and Dedner et al. \cite{Dedneretal}. Since we already use a staggered mesh for 
the semi-implicit discretization of the pressure subsystem, it is very natural to employ the strategy of Balsara and Spicer \cite{BalsaraSpicer1999,Balsara2004}, 
which also adopts a staggered mesh for the time evolution of the magnetic field. In this paper, we properly extend this technique to deal also with the resistive 
terms. For that purpose, it has to be noted that with $\nabla \cdot \mathbf{B}=0$ the resistive term can be rewritten in terms of a double curl operator as  
$\eta \nabla \cdot \left( \nabla \mathbf{B} - \nabla \mathbf{B}^T \right)  = - \eta \nabla \times \nabla \times \mathbf{B}$ and therefore 
the induction equation for the magnetic field reads 
\begin{equation}
  \frac{\partial \mathbf{B}}{\partial t} + \nabla \times \mathbf{E} = 0, 
\label{eqn.induction} 
\end{equation} 
with the electric field vector given by 
\begin{equation}
 \mathbf{E} = - \mathbf{v} \times \mathbf{B} + \eta \nabla \times \mathbf{B}, 
\label{eqn.Evector} 
\end{equation}
which reduces to the standard expression $\mathbf{E} = - \mathbf{v} \times \mathbf{B}$ for the ideal MHD equations ($\eta=0$).   
Again we split the MHD system into a first subsystem that contains the convective and the viscous terms that will both be discretized \textit{explicitly}, 
while the second one is again the pure pressure subsystem that will be discretized \textit{implicitly}, i.e. we write 
\begin{equation} 
   \frac{\partial \Q}{\partial t} + \nabla \cdot \left( \mathbf{F}^c - \mathbf{F}^v \right) + \nabla \cdot \mathbf{F}^p = 0, 
	\label{eqn.split2d} 
\end{equation} 
with 
\begin{equation} 
	\mathbf{F}^c = (\mathbf{f}^c, \mathbf{g}^c) = \left( \begin{array}{c} \rho \mathbf{v} \\ \rho \mathbf{v} \otimes \mathbf{v} +  m \mathbf{I} - \frac{1}{4\pi} \mathbf{B} \otimes \mathbf{B} \\  \mathbf{v}^T \left( \rho k + 2 m \right)  - \frac{1}{4\pi} \mathbf{v}^T \mathbf{B} \otimes \mathbf{B} \\ \mathbf{B} \otimes \mathbf{v} - \mathbf{v} \otimes \mathbf{B} \end{array}  \right), 
	\quad 
	\mathbf{F}^p = ( \mathbf{f}^p, \mathbf{g}^p ) = \left(\begin{array}{c} 0 \\ p  \mathbf{I} \\ (\rho \mathbf{v}) h \\ \mathbf{0} \end{array} \right) . 
\end{equation} 

\subsection{Semi-implicit discretization} 
\label{sec.si2d} 

The computational domain $\Omega$ is discretized by the control volumes of a primary grid denoted by $\Omega_{i,j}=[x_{i-\halb},x_{i+\halb}] \times [y_{j-\halb},y_{j+\halb}]$. 
To ease notation, in the following we suppose an equidistant mesh spacing of size $\Delta x$ and $\Delta y$ in $x$ and $y$ direction and the corresponding number of 
cells is denoted by $N_x$ and $N_y$, respectively. The edge-based staggered dual control volumes in $x$ direction are denoted by 
$\Omega_{i+\halb,i} = [x_i,x_{i+1}] \times [y_{j-\halb},y_{j+\halb}]$, while the control volumes of the staggered dual grid in $y$ direction are
$\Omega_{i,j+\halb} = [x_{i-\halb},x_{i+\halb}] \times [y_j,y_{j+1}]$, i.e. overall the method uses a set of three overlapping grids, each of which entirely
covers the domain $\Omega$. The averaging operators from the main grid to the dual grids read 
\begin{equation}
   \Q_{i+\halb,j}^n = \halb \left( \Q_{i,j}^n + \Q_{i+1,j}^n \right), \qquad 
   \Q_{i,j+\halb}^n = \halb \left( \Q_{i,j}^n + \Q_{i,j+1}^n \right), 
\end{equation} 
while the averaging from the two dual grids to the main grid is given by 
\begin{equation} 
   \Q_{i,j}^n = \halb \left( \Q_{i-\halb,j}^n + \Q_{i+\halb,j}^n \right), \qquad 
   \Q_{i,j}^n = \halb \left( \Q_{i,j-\halb}^n + \Q_{i,j+\halb}^n \right). 
\end{equation} 
\subsubsection{Convective and viscous subsystem.} 
The viscous and convective subsystem is discretized with an explicit finite volume scheme of the type 
\begin{equation}
   \Q_{i,j}^* = \Q_{i,j}^n - \frac{\Delta t}{\Delta x} \left( \mathbf{f}_{i+\halb,j} - \mathbf{f}_{i-\halb,j} \right) - \frac{\Delta t}{\Delta y} \left( \mathbf{g}_{i,j+\halb} - \mathbf{g}_{i,j-\halb} \right),     
	\label{eqn.expfv2d} 
\end{equation} 
where the numerical fluxes at the element interfaces contain both the nonlinear convective as well as the viscous terms and therefore read
\begin{eqnarray}
  \mathbf{f}_{i+\halb,j} &=& \halb \left( \mathbf{f}^c( \mathbf{Q}_{i+\halb,j}^-) + \mathbf{f}^c( \mathbf{Q}_{i+\halb,j}^+) \right) - 
	   \halb s_{\max}^x \left( \mathbf{Q}_{i+\halb,j}^+  -  \mathbf{Q}_{i+\halb,j}^- \right) -  \nonumber \\
			&&  \halb \left( \mathbf{f}^v(\mathbf{V}_{i+\halb,j+\halb}^n, \nabla \mathbf{V}_{i+\halb,j+\halb}^n) + \mathbf{f}^v(\mathbf{V}_{i+\halb,j-\halb}^n, \nabla \mathbf{V}_{i+\halb,j-\halb}^n) \right), 
\end{eqnarray} 
with the boundary extrapolated values $\mathbf{Q}_{i+\halb,j}^\pm$ and the maximum signal speed in $x$ direction $s_{\max}^x$ computed as in the one-dimensional case. The expression for the numerical 
flux $\mathbf{g}_{i,j+\halb}$ is obviously very similar to the one for $\mathbf{f}_{i+\halb,j}$, hence it is not necessary to report it here.  
The second order MUSCL-Hancock scheme in two space dimensions is a straight-forward extension of the one-dimensional case shown previously and is well known, so we can omit the details. 
For the viscous flux we define the \textit{corner variables} 
\begin{equation}
   \mathbf{V}_{i+\halb,j+\halb}^n = \frac{1}{4} \left( \mathbf{V}_{i,j}^n +  \mathbf{V}_{i+1,j}^n +  \mathbf{V}_{i,j+1}^n +  \mathbf{V}_{i+1,j+1}^n \right), 
\label{eqn.corner.value}	
\end{equation}
and the \textit{corner gradients} of the vector of primitive variables given by 
\begin{equation} 
    \nabla \mathbf{V}_{i+\halb,j+\halb}^n = \halb \left(   
     \frac{ \mathbf{V}_{i+1,j+1}^n - \mathbf{V}_{i,j+1}^n}{\Delta x}  + \frac{ \mathbf{V}_{i+1,j}^n - \mathbf{V}_{i,j}^n}{\Delta x} , 
     \frac{ \mathbf{V}_{i+1,j+1}^n - \mathbf{V}_{i+1,j}^n}{\Delta y}  + \frac{ \mathbf{V}_{i,j+1}^n - \mathbf{V}_{i,j}^n}{\Delta y}  
		                 \right). 
\label{eqn.corner.gradient}
\end{equation}  
\subsubsection{Divergence-free evolution of the magnetic field.}
In multiple space dimensions, it is of fundamental importance to evolve the magnetic field in a consistent manner that respects the divergence-free condition $\nabla \cdot \mathbf{B} = 0$ exactly
also on the discrete level. For this purpose, we follow \cite{BalsaraSpicer1999,Balsara2004} and introduce the magnetic field components on the staggered mesh as $(B_x)_{i+\halb,j}^n$
and $(B_y)_{i,j+\halb}^n$. The normal magnetic field components can then be evolved in time by a discrete form of the induction equation \eqref{eqn.induction} as follows
\begin{equation}
  (B_x)_{i+\halb,j}^{n+1} = (B_x)_{i+\halb,j}^n - \frac{\Delta t}{\Delta y} \left( E^z_{i+\halb,j+\halb} - E^z_{i+\halb,j-\halb} \right), 
	\label{eqn.bx} 
\end{equation}
\begin{equation}	
  (B_y)_{i,j+\halb}^{n+1} = (B_y)_{i,j+\halb}^n + \frac{\Delta t}{\Delta x} \left( E^z_{i+\halb,j+\halb} - E^z_{i-\halb,j+\halb} \right), 	
	\label{eqn.by} 
\end{equation}
with the electric field component in $z$ direction given by a multi-dimensional Riemann solver (see e.g. \cite{balsarahlle2d,balsarahllc2d,balsarahlle3d,BalsaraMultiDRS,MUSIC1,MUSIC2}) as 
\begin{eqnarray}
  E^z_{i+\halb,j+\halb} &=& \phantom{-} \halb v_{i+\halb,j+\halb}^n \left( (B_x)^n_{i+\halb,j}   + (B_x)^n_{i+\halb,j+1} \right) 
	                        - \halb s_{\max}^y                        \left( (B_x)^n_{i+\halb,j+1} - (B_x)^n_{i+\halb,j}   \right) \nonumber \\ 
									      &&  - \halb u_{i+\halb,j+\halb}^n 	        \left( (B_y)^n_{i,j+\halb}   + (B_y)^n_{i+1,j+\halb} \right) 
													  + \halb s_{\max}^x                      \left( (B_y)^n_{i+1,j+\halb} - (B_y)^n_{i,j+\halb}   \right) \nonumber \\
												&&	+ \eta \left( \partial_x (B_y)_{i+\halb,j+\halb}^n - \partial_y (B_x)_{i+\halb,j+\halb}^n \right). 
 	\label{eqn.ez} 
\end{eqnarray} 
Note that in \eqref{eqn.ez} the last line accounts for the resistive term and is an approximation to the $z$ component of the curl of $\mathbf{B}$ using
the corner gradients computed in \eqref{eqn.corner.gradient}. The velocity vector in the corner has already been computed via 
\eqref{eqn.corner.value}. 
It is easy to check that the scheme \eqref{eqn.bx}-\eqref{eqn.by} is exactly divergence-free in the discrete sense 
\begin{equation}
   \frac{(B_x)_{i+\halb,j}^{n+1}-(B_x)_{i-\halb,j}^{n+1}}{\Delta x} + \frac{(B_y)_{i,j+\halb}^{n+1}-(B_y)_{i,j-\halb}^{n+1}}{\Delta y} = 0, 
	\label{eqn.disc.divb} 
\end{equation} 
if the magnetic field was discretely divergence-free at the initial time $t=0$. Note that in 2D it is sufficient to take $(B_z)_{i,j}^{n+1}=(B_z)_{i,j}^*$ from \eqref{eqn.expfv2d}. 
After the update of the staggered magnetic fields $B_x$ and $B_y$ via \eqref{eqn.bx}-\eqref{eqn.by}, the cell-centered magnetic field vector 
$\mathbf{B}_{i,j}^{n+1}$ is obtained by averaging the staggered quantities back from the dual grid to the main grid. It has to be stressed that in the multidimensional case 
in general $\mathbf{B}_{i,j}^{n+1} \neq \mathbf{B}_{i,j}^{*}$, i.e. the cell-centered quantity $\mathbf{B}_{i,j}^{*}$ obtained from \eqref{eqn.expfv2d} is 
only an auxiliary quantity that is \textit{overwritten} by the averages onto the main grid of the consistently evolved magnetic field components $(B_x)_{i+\halb,j}^{n+1}$ 
and $(B_y)_{i,j+\halb}^{n+1} $, which are the main quantities that represent the discrete magnetic field in our scheme. The cell-centered magnetic field is needed in order
to compute the energy density of the magnetic field $m_{i,j}^{n+1}$ needed later in the pressure subsystem. 

\subsubsection{Pressure subsystem.} 
In two space dimensions the discrete momentum equations read 
\begin{equation}
\label{eqn.rhou2d} 
 (\rho u)_{i+\halb,j}^{n+1} = (\rho u)_{i+\halb,j}^{*} - \frac{\Delta t}{\Delta x} \left( p_{i+1,j}^{n+1} - p_{i,j}^{n+1} \right),  \quad 
 (\rho v)_{i,j+\halb}^{n+1} = (\rho v)_{i,j+\halb}^{*} - \frac{\Delta t}{\Delta y} \left( p_{i,j+1}^{n+1} - p_{i,j}^{n+1} \right),  
\end{equation} 
where pressure is taken \textit{implicitly}, while all nonlinear convective and viscous terms have already been discretized \textit{explicitly} via the operators $(\rho u)_{i+\halb,j}^{*} $ and 
$(\rho v)_{i,j+\halb}^{*}$ given in \eqref{eqn.expfv2d}. A preliminary form of the discrete total energy equation reads 
\begin{eqnarray}
\label{eqn.rhoE2d.prelim} 
  \rho e\left( p_{i,j}^{n+1}, \rho_{i,j}^{n+1} \right) + m_{i,j}^{n+1} + (\tilde{\rho k})^{n+1}_{i,j}  = (\rho E)_{i,j}^* 
\nonumber \\      
- \frac{\Delta t}{\Delta x} \left(  \tilde{h}_{i+\halb,j}^{n+1} (\rho u)_{i+\halb,j}^{n+1} - \tilde{h}_{i-\halb,j}^{n+1} (\rho u)_{i-\halb,j}^{n+1} \right)  
- \frac{\Delta t}{\Delta y} \left(  \tilde{h}_{i,j+\halb}^{n+1} (\rho v)_{i,j+\halb}^{n+1} - \tilde{h}_{i,j-\halb}^{n+1} (\rho v)_{i,j-\halb}^{n+1} \right). \nonumber \\  
\end{eqnarray}
Here, we have used again the abbreviation $\rho e\left( p_{i,j}^{n+1}, \rho_{i,j}^{n+1} \right) = \rho_{i,j}^{n+1} e\left( p_{i,j}^{n+1}, \rho_{i,j}^{n+1} \right)$.   
Inserting the discrete momentum equations \eqref{eqn.rhou2d} into the discrete energy equation \eqref{eqn.rhoE2d.prelim} and making tilde symbols explicit via the simple 
Picard iteration, as in the one-dimensional case, leads to the following discrete wave equation for the unknown pressure: 
\begin{eqnarray}
\label{eqn.p2d} 
  \rho_{i,j}^{n+1} e\left( p_{i,j}^{n+1,r+1}, \rho_{i,j}^{n+1} \right) & & \nonumber \\
- \frac{\Delta t^2}{\Delta x^2} \left(  {{h}_{i+\halb,j}^{n+1,r}} \left( p_{i+1,j}^{n+1,r+1}-p_{i,j}^{n+1,r+1} \right) 
                              - {{h}_{i-\halb,j}^{n+1,r}} \left( p_{i,j}^{n+1,r+1}-p_{i-1,j}^{n+1,r+1} \right) \right) & & \nonumber \\  
- \frac{\Delta t^2}{\Delta y^2}\left(  {{h}_{i,j+\halb}^{n+1,r}} \left( p_{i,j+1}^{n+1,r+1}-p_{i,j}^{n+1,r+1} \right) 
                              - {{h}_{i,j-\halb}^{n+1,r}} \left( p_{i,j}^{n+1,r+1}-p_{i,j-1}^{n+1,r+1} \right) \right)  
 & = & b_{i,j}^r, 
\end{eqnarray} 
with the known right hand side 
\begin{eqnarray}
 b_{i,j}^r  =   (\rho E)_{i,j}^* - m_{i,j}^{n+1} - (\rho k)^{n+1,r}_{i,j}   
 \nonumber \\
   - \frac{\Delta t}{\Delta x} \left( {h}_{i+\halb,j}^{n+1,r} (\rho u)_{i+\halb,j}^{*} - {h}_{i-\halb,j}^{n+1,r} (\rho u)_{i-\halb,j}^{*} \right)
   - \frac{\Delta t}{\Delta y} \left( {h}_{i,j+\halb}^{n+1,r} (\rho v)_{i,j+\halb}^{*} - {h}_{i,j-\halb}^{n+1,r} (\rho v)_{i,j-\halb}^{*} \right). 
\label{eqn.rhs.2d} 
\end{eqnarray} 
We stress that the density $\rho_{i,j}^{n+1} = \rho_{i,j}^{*} $ is already known from \eqref{eqn.expfv2d}, and also the energy of the magnetic field $m_{i,j}^{n+1}$ is 
already known after averaging the staggered normal magnetic field components that have been evolved via \eqref{eqn.bx} and \eqref{eqn.by} onto the main grid. 
The system for the pressure \eqref{eqn.p2d} is again a mildly nonlinear system of the form \eqref{eqn.nonlinear} with a linear part that is symmetric and as 
least positive semi-definite. Hence, with the usual assumptions on the nonlinearity detailed in \cite{CasulliZanolli2012}, it can be again efficiently solved with 
the nested Newton method of Casulli and Zanolli \cite{CasulliZanolli2010,CasulliZanolli2012}. Note that in the incompressible limit $M \to 0$, following the 
asymptotic analysis performed in \cite{KlaMaj,KlaMaj82,Klein2001,Munz2003,MunzDumbserRoller}, the pressure tends to a constant and the contribution of the kinetic energy $\rho k$ 
can be neglected w.r.t. $\rho e$. Therefore, in the incompressible limit the system \eqref{eqn.p2d} tends to the usual pressure Poisson equation of incompressible 
flow solvers. 
In each Picard iteration, after the solution of the pressure system \eqref{eqn.p2d} the enthalpies at the interfaces can be updated and the momentum is updated by 
\begin{eqnarray}
\label{eqn.rhou2d.pic} 
 (\rho u)_{i+\halb,j}^{n+1,r+1} &=& (\rho u)_{i+\halb,j}^{*} - \frac{\Delta t}{\Delta x} \left( p_{i+1,j}^{n+1,r+1} - p_{i,j}^{n+1,r+1} \right),  \\
 (\rho v)_{i,j+\halb}^{n+1,r+1} &=& (\rho v)_{i,j+\halb}^{*} - \frac{\Delta t}{\Delta y} \left( p_{i,j+1}^{n+1,r+1} - p_{i,j}^{n+1,r+1} \right),  
\end{eqnarray} 
from which $(\rho k)_{i,j}^{n+1,r+1}$ can be computed after averaging onto the main grid.  
At the end of the Picard iterations, the total energy is updated as 
\begin{eqnarray}
\label{eqn.rhoE2d} 
 (\rho E)_{i,j}^{n+1} &=& (\rho E)_{i,j}^* 
- \frac{\Delta t}{\Delta x} \left( {h}_{i+\halb,j}^{n+1} (\rho u)_{i+\halb,j}^{n+1} - {h}_{i-\halb,j}^{n+1} (\rho u)_{i-\halb,j}^{n+1} \right) \nonumber \\
& & \phantom{(\rho E)_{i,j}^*} 
- \frac{\Delta t}{\Delta y} \left( {h}_{i,j+\halb}^{n+1} (\rho v)_{i,j+\halb}^{n+1} - {h}_{i,j-\halb}^{n+1} (\rho v)_{i,j-\halb}^{n+1} \right),   
\end{eqnarray} 
while the final momentum is averaged back onto the main grid. This completes the description of our new divergence-free semi-implicit algorithm 
for the VRMHD equations in the multi-dimensional case. 

\section{Numerical results in 2D}
\label{sec.results2d} 

In all the following numerical test problems, the ideal gas equation of state is used, in order to make the results comparable with existing 
data in the literature. For applications with general EOS, see \cite{DumbserCasulli2016}. 
If not specified otherwise, the ratio of specific heats is chosen as $\gamma=1.4$ in all the following test cases. 
The CPU timings reported in this section were obtained on a workstation using one single core of an Intel i7-2600 CPU with 3.4 GHz clock 
speed and 12 GB of RAM. In order to allow a better quantitative comparison with other schemes, we report the average CPU time that was
needed to carry out one time step for one control volume, i.e. dividing the total wall clock time needed by the simulation by the number of
time steps and the number of control volumes. The inverse of this number corresponds to the number of zones which the scheme is able to 
update within one second of wallclock time on one CPU core. In the 2D simulations the time step is computed according to 
\begin{equation}
 \Delta t = \CFL \frac{1}{ \frac{\max |\lambda^c_x|}{\Delta x} + \frac{ \max |\lambda^c_y|}{\Delta y} + 2 \left( \frac{4}{3} \frac{\mu}{\rho} + \frac{\lambda}{c_v \rho} + \eta \right) \left( \frac{1}{\Delta x^2} + \frac{1}{\Delta y^2} \right) },
\end{equation} 
with the Courant number $\CFL < 1$ and the 'convective' eigenvalues $\lambda^c_x$ and $\lambda^c_y$ in $x$ and
$y$ direction, respectively. If not specified otherwise, we set $\CFL=0.9$ in all test problems presented in this section. Furthermore, for all test 
cases we have explicitly verified that up to machine precision the magnetic field is divergence-free and mass, momentum and energy are conserved. 

\subsection{Low Mach number magnetic field loop advection} 
\label{sec.fieldloop}

\begin{figure}[!htbp]
\begin{center}
\begin{tabular}{cc} 
\includegraphics[width=0.45\textwidth]{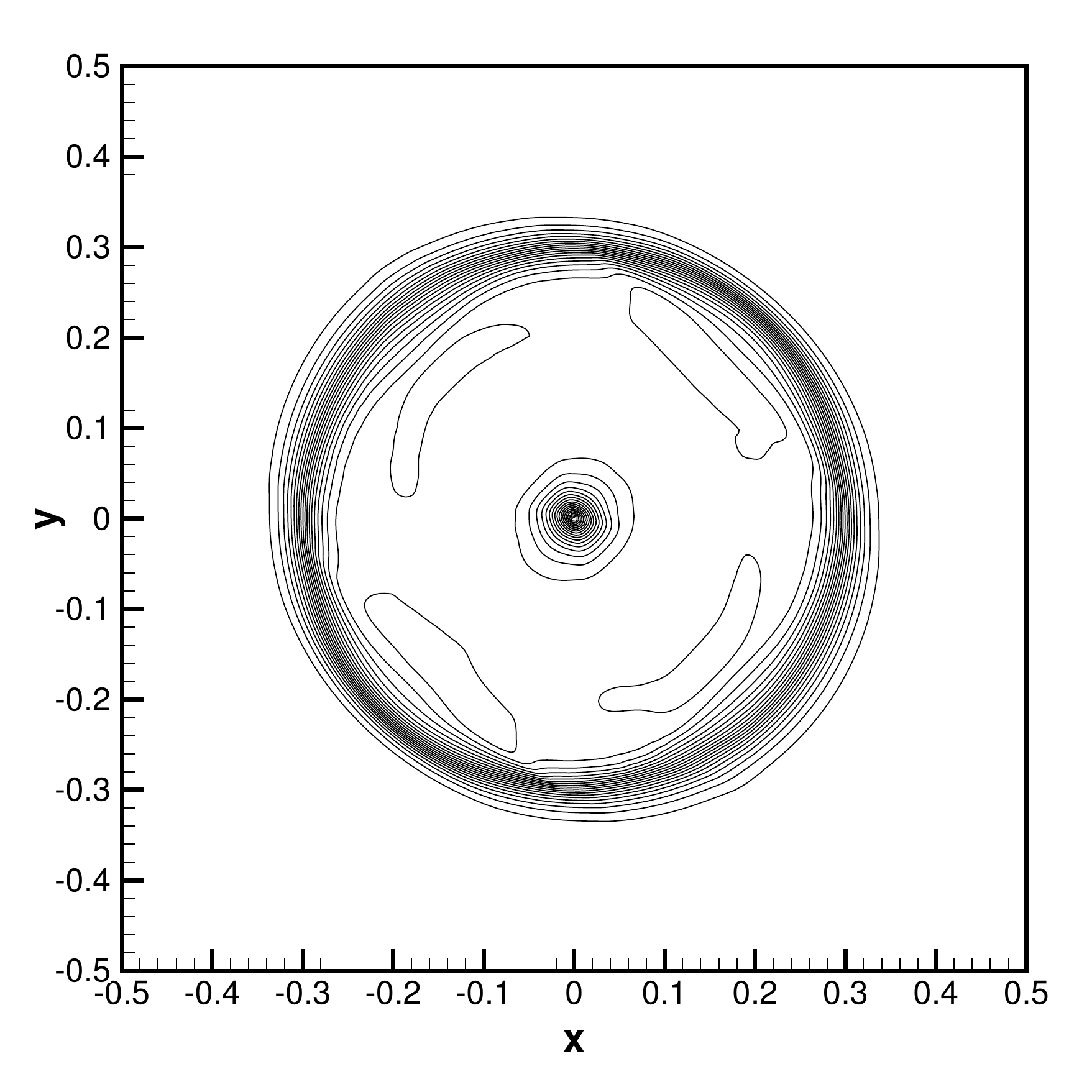}   & 
\includegraphics[width=0.45\textwidth]{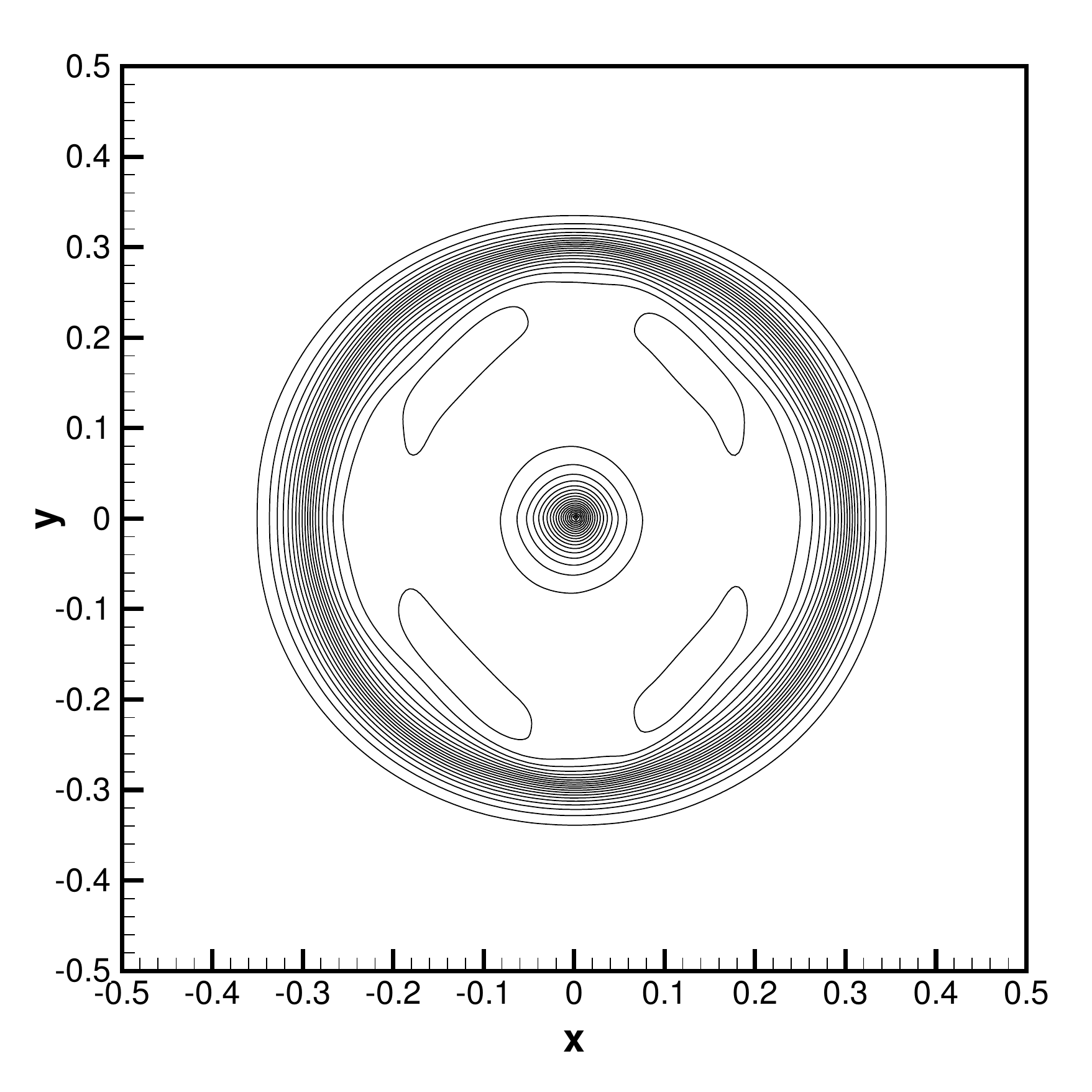}     
\end{tabular} 
\caption{Numerical solution at time $t=1.0$ obtained for the low Mach number MHD field loop advection problem with the divergence-free SIFV method (left) and with 
a divergence-free explicit second order Godunov-type scheme (right). 20 equidistant contour lines of the magnetic field strength in the interval $[10^{-5},10^{-3}]$ are shown. For this test
case, the density-based explicit scheme was more than a factor of 50 slower than the new pressure-based semi-implicit method. } 
\label{fig.loop}
\end{center}
\end{figure}

Here we solve the magnetic field loop advection problem proposed by Gardiner and Stone in \cite{GardinerStone}. However, in order to make it more difficult 
and in order to show the performance of our new divergence-free semi-implicit finite volume scheme, we run the test case at \textit{low Mach number}. 
The setup of the test problem is described in the following. The computational domain is $\Omega=[-1,1] \times [-\halb,\halb]$ with four periodic boundary 
conditions everywhere. The initial density is set to $\rho=1$, the initial velocity field is $\mathbf{v}=(2,1,0)$, the pressure is $p=10^5$ and the 
initial magnetic field is prescribed by the magnetic vector potential 
\begin{equation}
  A = \left\{ \begin{array}{lll} A_0 (R-r) & \textnormal{ if } & r \leq R, \\ 
	                                 0       & \textnormal{ if } & r > R,
	\end{array} \right. 
\end{equation} 
with $A_0=10^{-3}$, $R=0.3$ and $r^2=x^2+y^2$. The Mach number of the flow is about $M=0.006$. We run the problem with the second order version of our new 
semi-implicit FV scheme and with a divergence-free second-order explicit MUSCL-type TVD finite volume scheme \cite{BalsaraSpicer1999,Balsara2004} until 
$t=1$ in order to complete one entire advection period. 
In both cases the domain $\Omega$ is discretized with $500 \times 250$ control volumes and the CFL number is set to $\CFL=0.8$.  
The computational results for both cases (explicit vs. semi-implicit) are depicted in Figure \ref{fig.loop} and are comparable with those obtained 
in the literature, see e.g. \cite{GardinerStone,balsarahlle2d,ADERdivB}, although the explicit scheme appears to be slightly more dissipative, probably 
due to the extremely large number of time steps needed to reach the final time. The explicit method needed a total wall-clock time of $78414$s, 
while our new semi-implicit FV scheme was able to complete the simulation in only $1356$s. This results in a speedup factor of \textbf{57} for 
the new semi-implicit scheme, which is a \textbf{clear advantage} for the new algorithm presented in this paper over existing schemes.  
For this simulation, the average computational cost of the SIFV scheme was $11.5 \mu$s per element and time step. The most 
expensive part here was the solution of the pressure system in the semi-implicit algorithm. For comparison, the average cost per element update for the 
explicit second order Godunov-type TVD scheme in this test was only $ 2\mu$s per element and time step. However, since the explicit scheme needs 
\textit{two orders of magnitude} more time steps compared to the semi-implicit scheme, the new SIFV method presented in this paper is still 
computationally much more efficient.

\subsection{Ideal MHD rotor problem}
\label{sec.rotor}

\begin{figure}[!htbp]
\begin{center}
\begin{tabular}{cc} 
\includegraphics[width=0.45\textwidth]{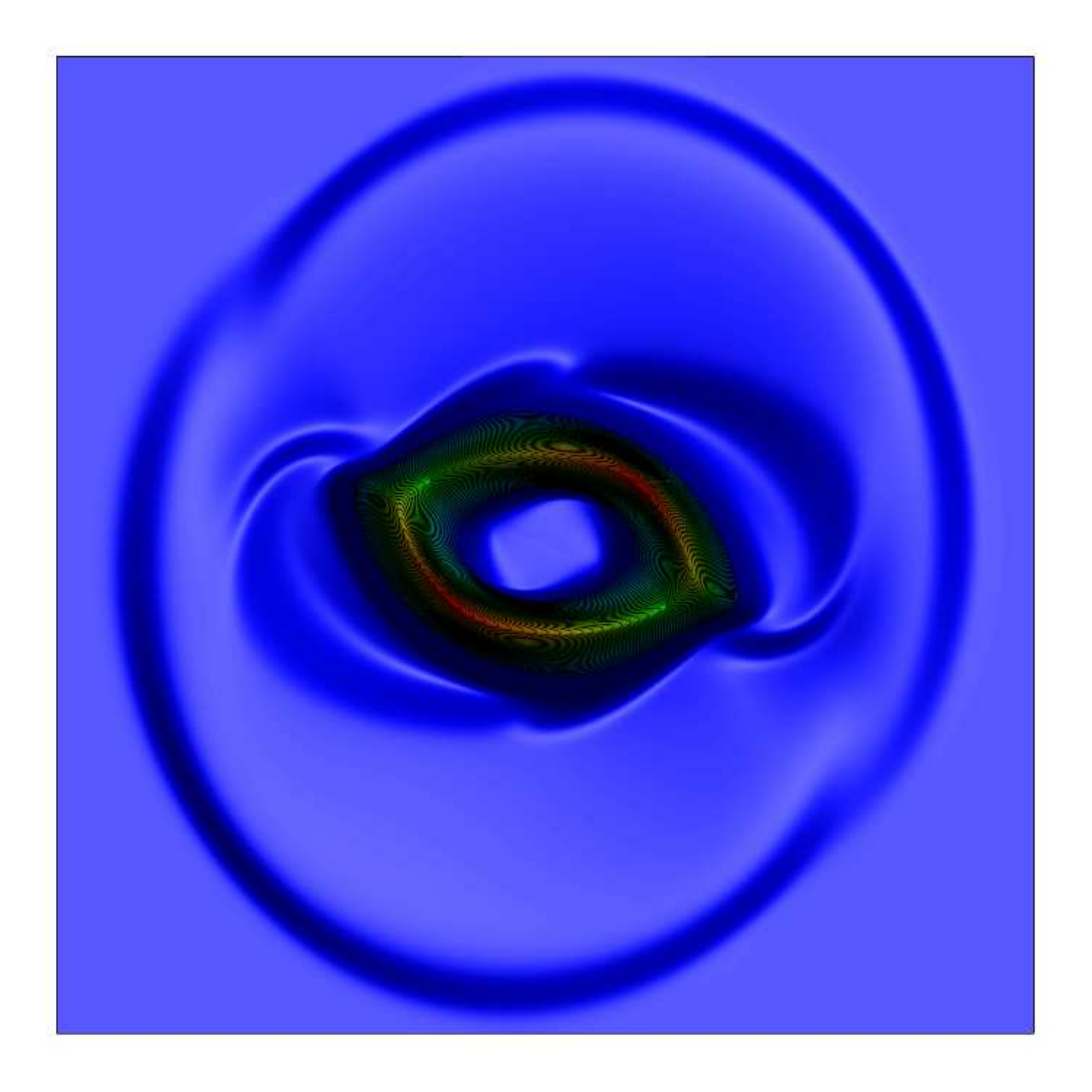} & 
\includegraphics[width=0.45\textwidth]{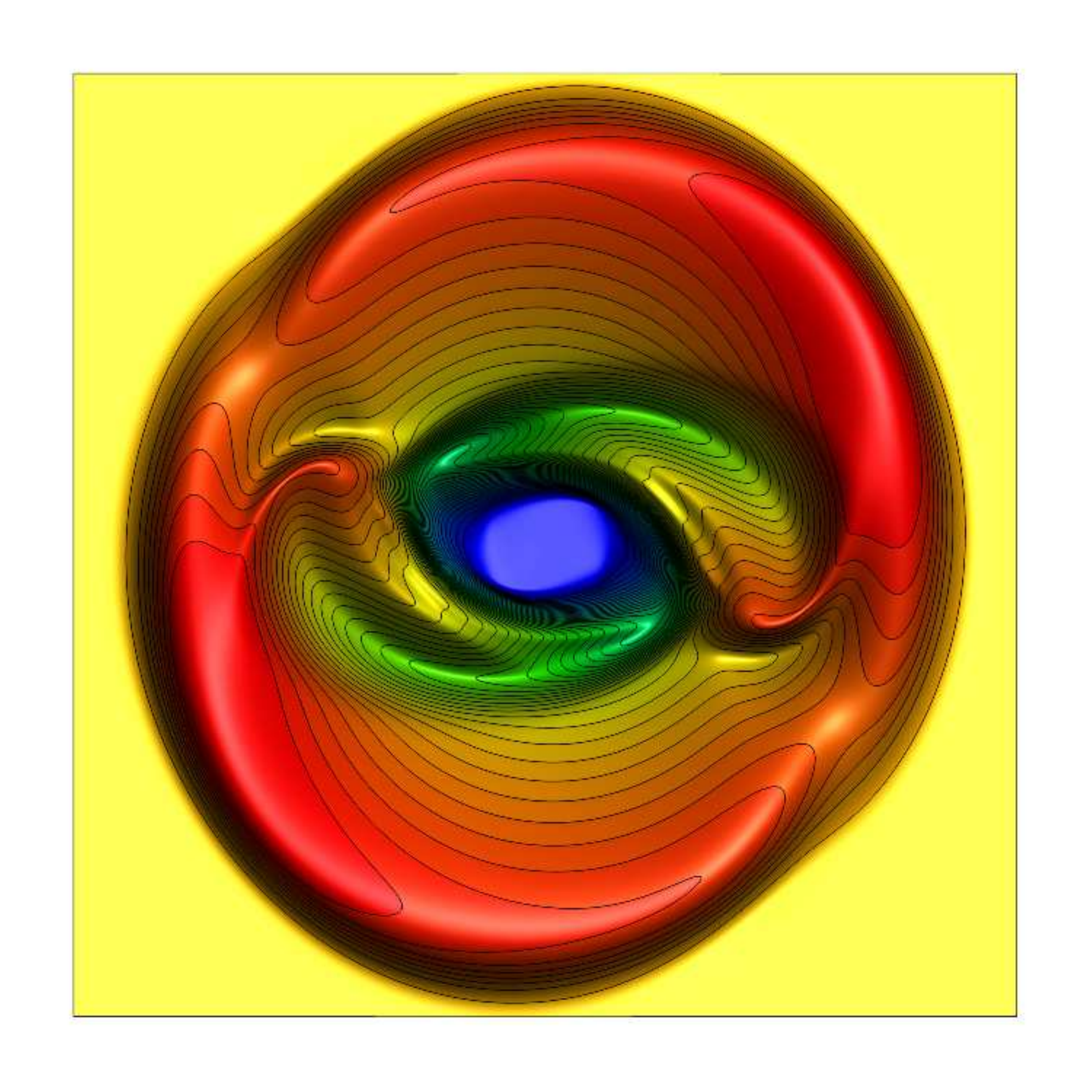}   \\ 
\includegraphics[width=0.45\textwidth]{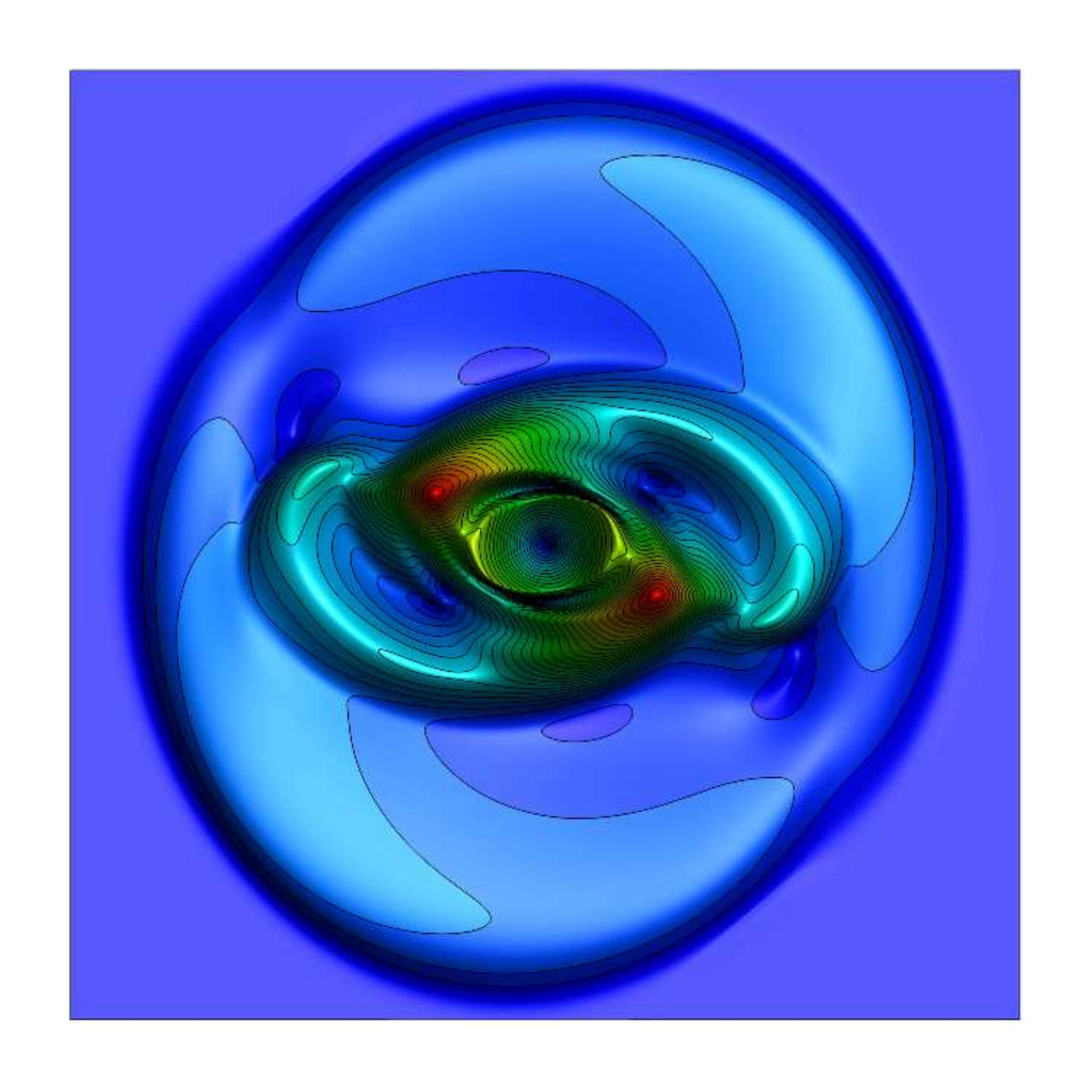} & 
\includegraphics[width=0.45\textwidth]{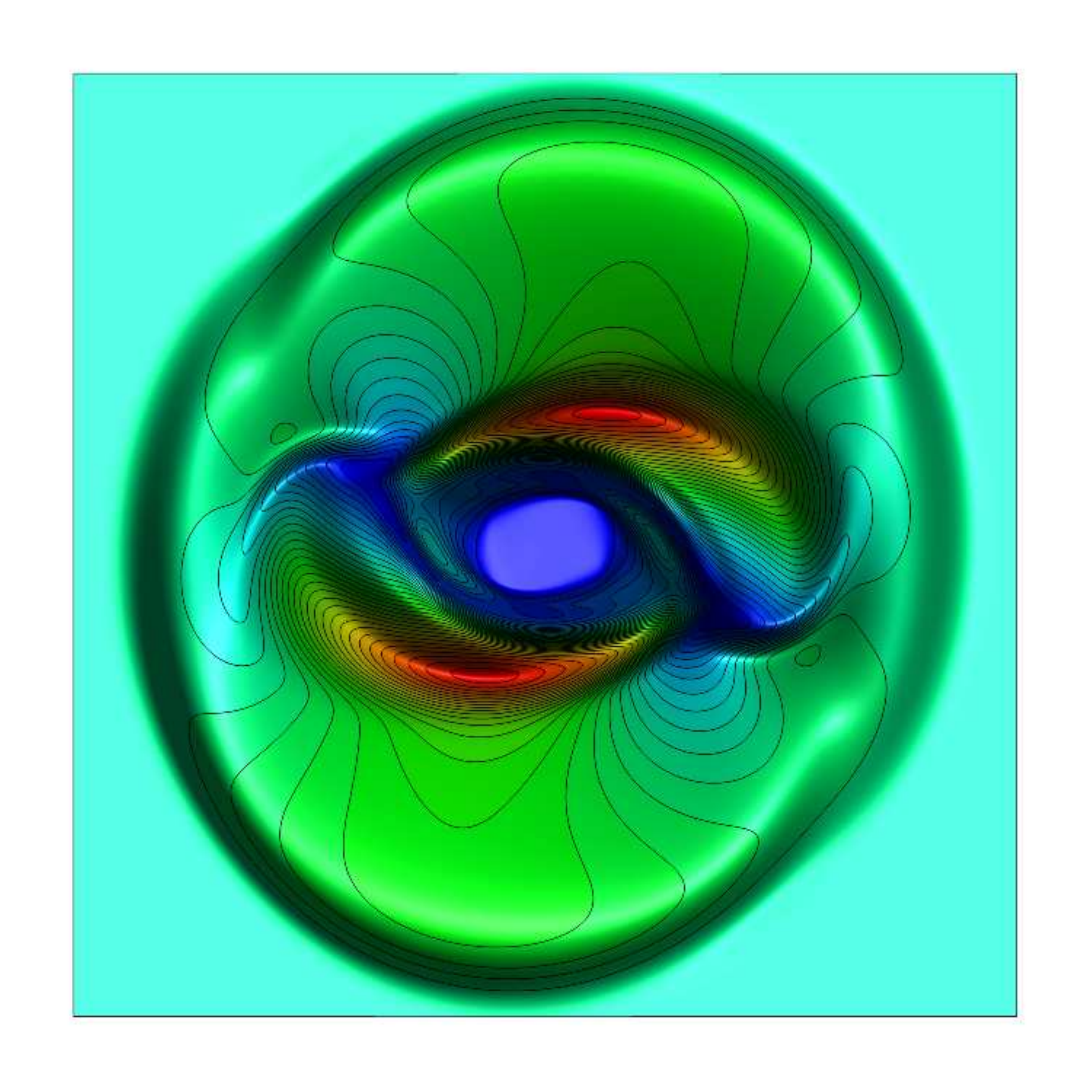}    
\end{tabular} 
\caption{Numerical solution obtained with the divergence-free semi-implicit finite volume method for the MHD rotor problem at time $t=0.25$. Contour lines of density (top left), pressure (top right), 
Mach number (bottom left) and magnetic pressure (bottom right). }  
\label{fig.rotor}
\end{center}
\end{figure}

The well-known MHD rotor problem of Balsara and Spicer \cite{BalsaraSpicer1999} has become a standard test bed for testing numerical methods for the ideal MHD equations. In this test a rotating high density fluid (the rotor) is embedded in a low density atmosphere at rest. Initially the pressure $p=1$ and the magnetic  field vector $\mathbf{B} = ( 2.5, 0, 0)^T$ are constant throughout the entire domain $\Omega = [-0.5,+0.5]^2$. 
The rotor produces torsional Alfv\'en waves which travel into the outer fluid at rest. The domain is discretized using a uniform Cartesian grid composed 
of $1000 \times 1000$ elements. For $0 \leq r \leq 0.1$, i.e. inside the rotor, the initial density is $\rho=10$, while it is set to $\rho=1$ outside. 
The velocity field inside the rotor is set to $\mathbf{v} = \boldsymbol{\omega} \times \mathbf{x}$ with $\boldsymbol{\omega}=(0,0,10)$, while $\mathbf{v}=(0,0,0)$
in the outer fluid. The computational results obtained with the new divergence-free semi-implicit finite volume scheme at time $t=0.25$ are shown in 
Fig. \ref{fig.rotor} for the fluid density, the pressure, the Mach number as well as the magnetic pressure. The results agree qualitatively well with those 
obtained by Balsara and Spicer in \cite{BalsaraSpicer1999} and other results reported elsewhere in the literature, see e.g. \cite{Dumbser2008,balsarahlle2d,AMR3DCL,ADERdivB,HPRmodelMHD}. 
The average computational cost of the SIFV scheme in this simulation was $3 \mu$s per element and time step. 

\subsection{Ideal MHD blast wave problem}
\label{sec.blast}

The MHD blast wave problem introduced in \cite{BalsaraSpicer1999} is a notoriously difficult test case. 
The initial data for density, velocity and magnetic field are constant throughout the domain and are set to $\rho=1$, $\mathbf{v}=(0,0,0)$ and $\mathbf{B}=(100,0,0)$. 
The pressure is initialized with $p=1000$ in an inner circular region $r<0.1$ and is set to $p=0.1$ outside, hence the pressure jumps over four orders of magnitude in this test problem. 
Furthermore, the fluid is highly magnetized due to the presence of a very strong magnetic field in the entire domain. The computational domain $\Omega=[-0.5,+0.5]^2$ 
is discretized with a uniform Cartesian grid using $1000 \times 1000$ pressure control volumes. The computational results obtained with our new divergence-free semi-implicit finite volume 
scheme at time $t=0.01$ are presented in Fig. \ref{fig.rotor} for the density, the pressure, the velocity magnitude and the magnetic pressure. The results agree qualitatively with those
obtained in the literature, see \cite{BalsaraSpicer1999,ADERdivB,HPRmodelMHD}.  
Also for this test problem the average computational cost of the SIFV scheme was $3 \mu$s per element and time step.

\begin{figure}[!htbp]
\begin{center}
\begin{tabular}{cc} 
\includegraphics[width=0.45\textwidth]{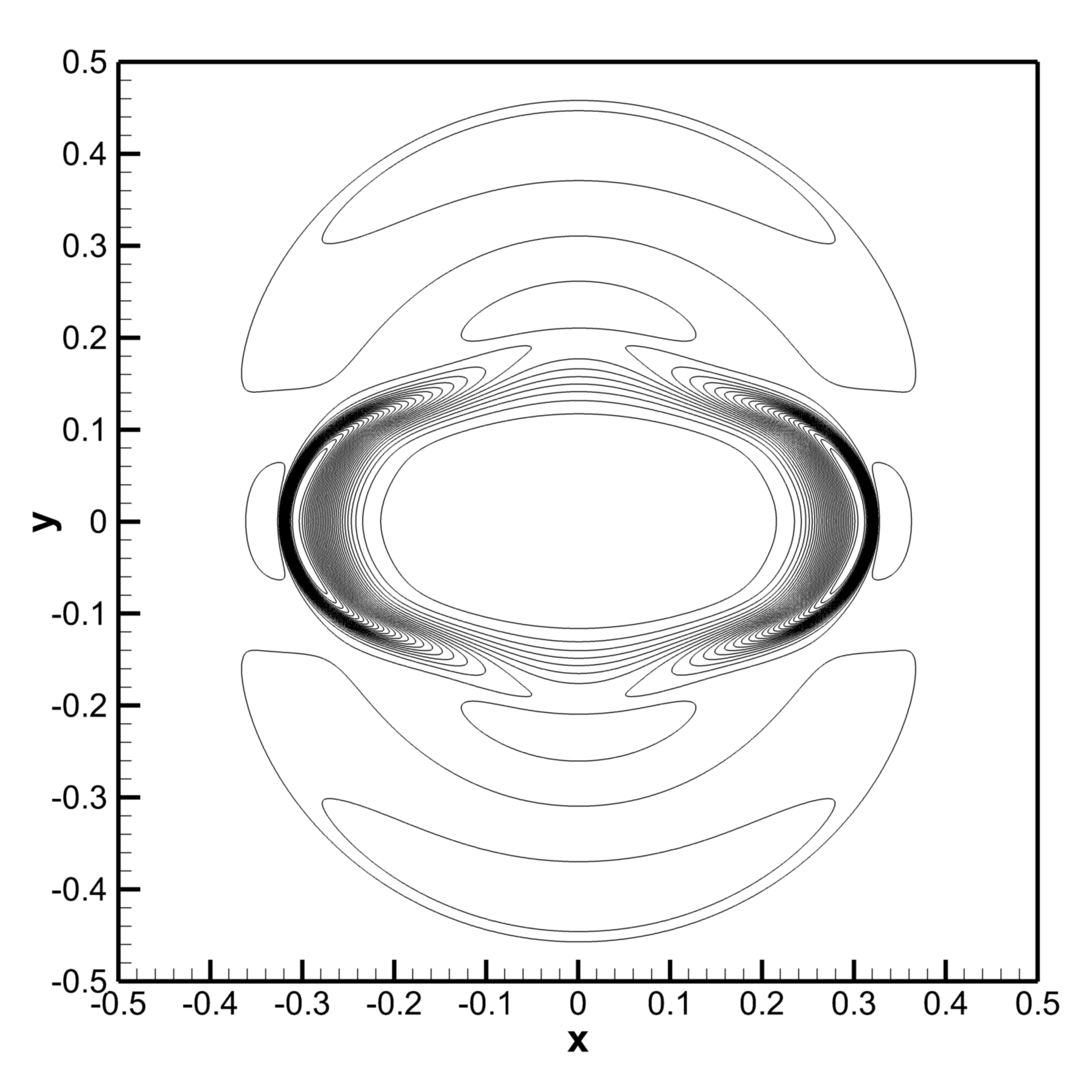} & 
\includegraphics[width=0.45\textwidth]{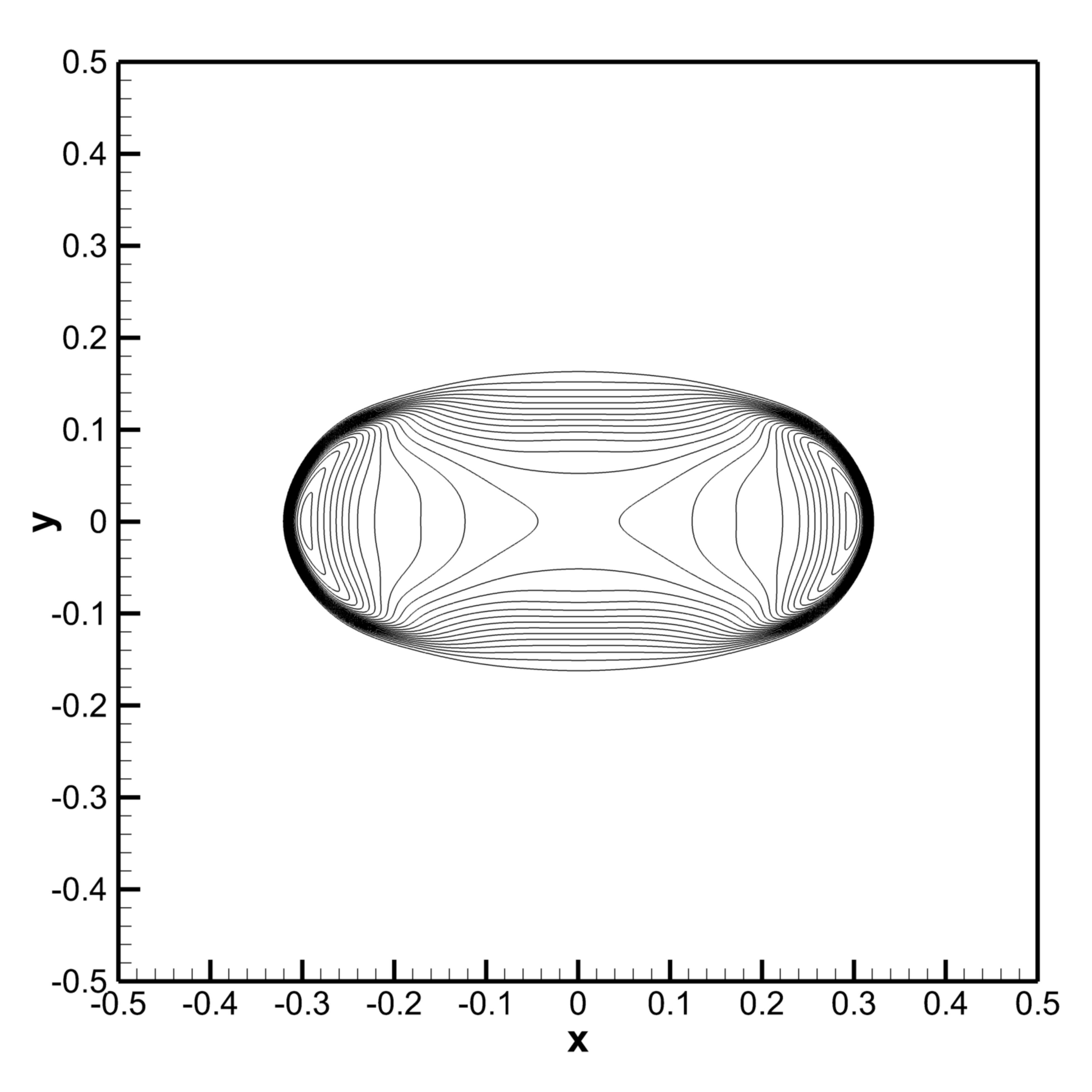}   \\ 
\includegraphics[width=0.45\textwidth]{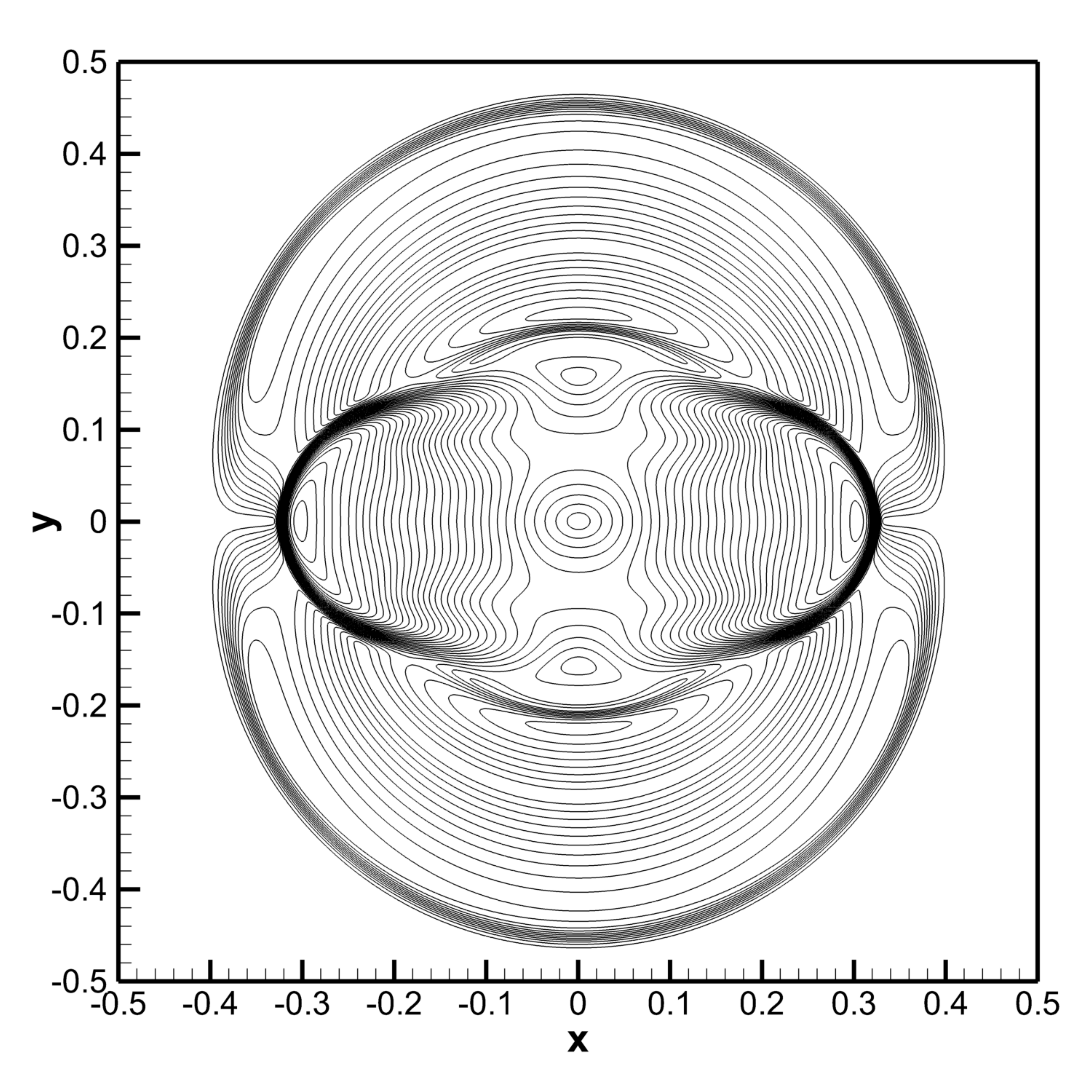} & 
\includegraphics[width=0.45\textwidth]{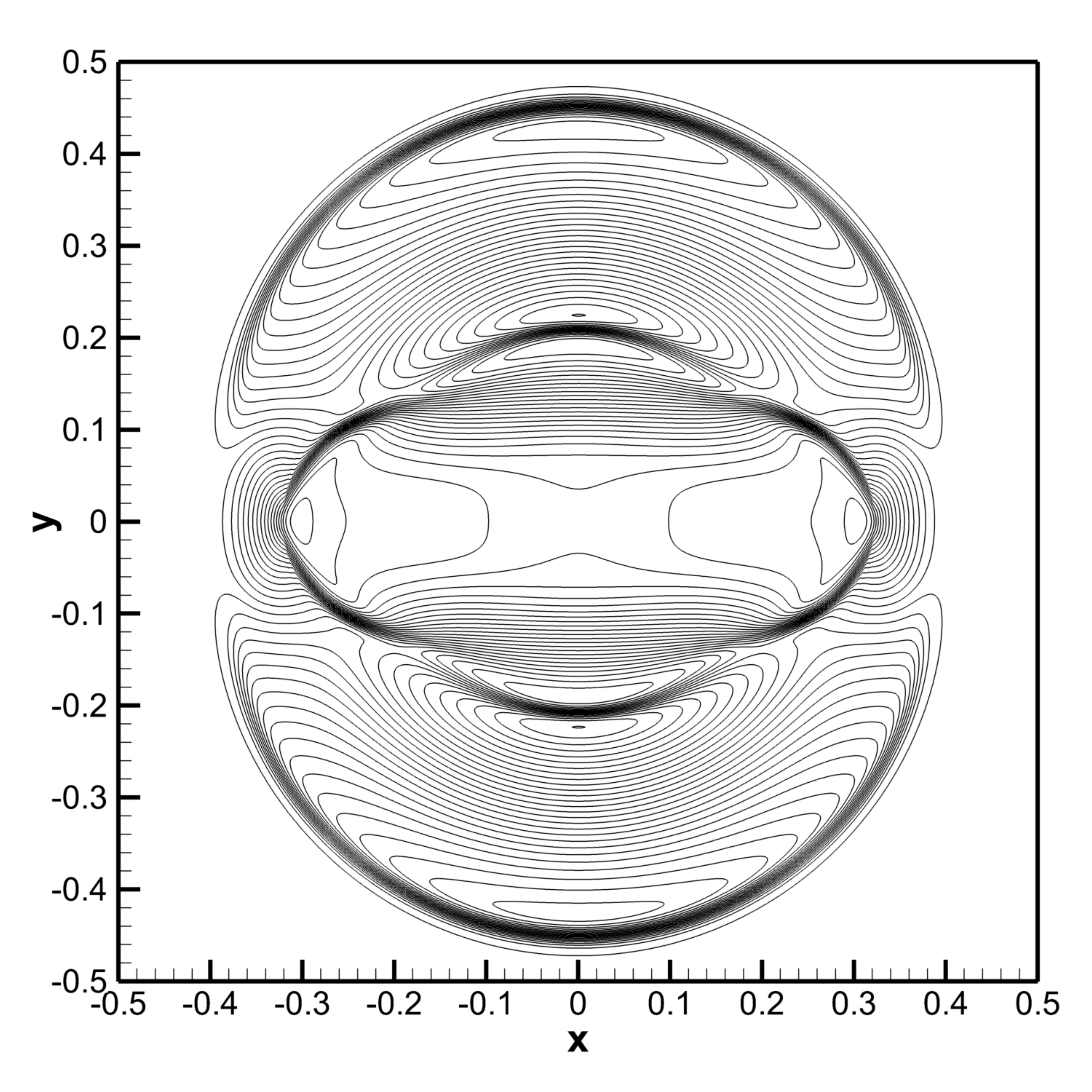}    
\end{tabular} 
\caption{Numerical solution obtained with the divergence-free semi-implicit finite volume method for the MHD blast wave problem at time $t=0.01$. Contour lines of density (top left), pressure (top right), velocity magnitude (bottom left) and magnetic pressure (bottom right). } 
\label{fig.blast}
\end{center}
\end{figure}

\subsection{Ideal MHD Orszag-Tang vortex}
\label{sec.ot}

\begin{figure}[!htbp]
\begin{center}
\begin{tabular}{cc} 
\includegraphics[width=0.45\textwidth]{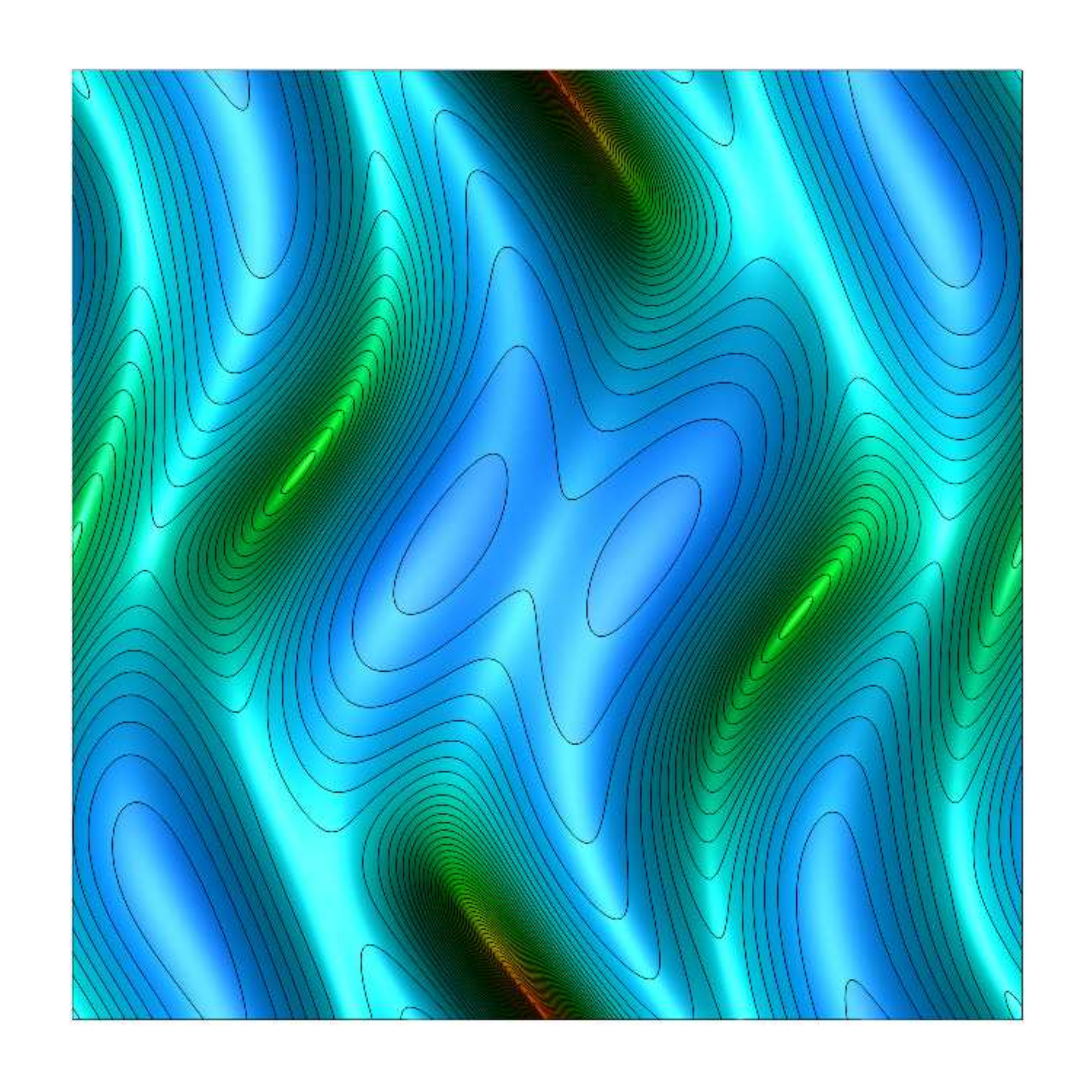} & 
\includegraphics[width=0.45\textwidth]{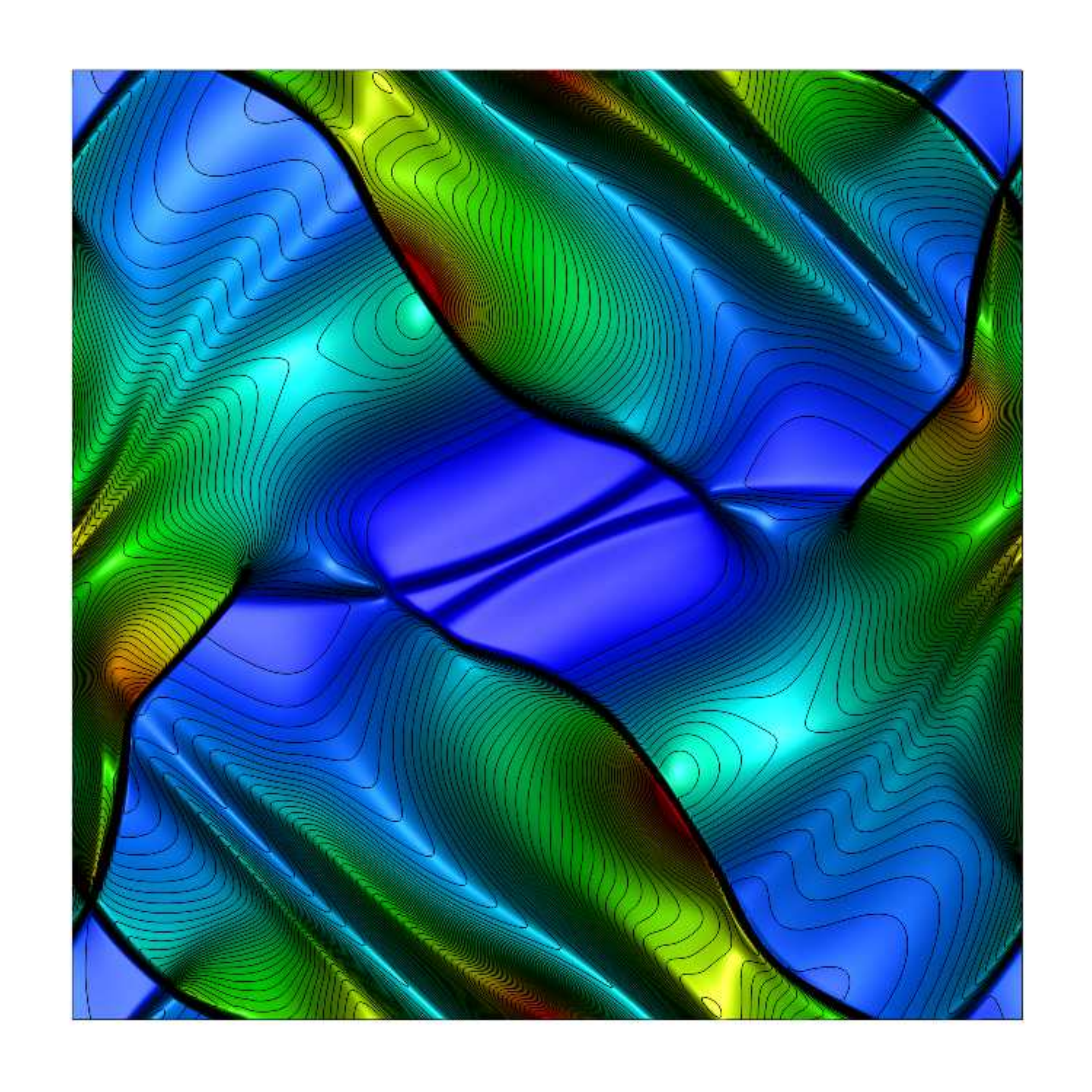}   \\ 
\includegraphics[width=0.45\textwidth]{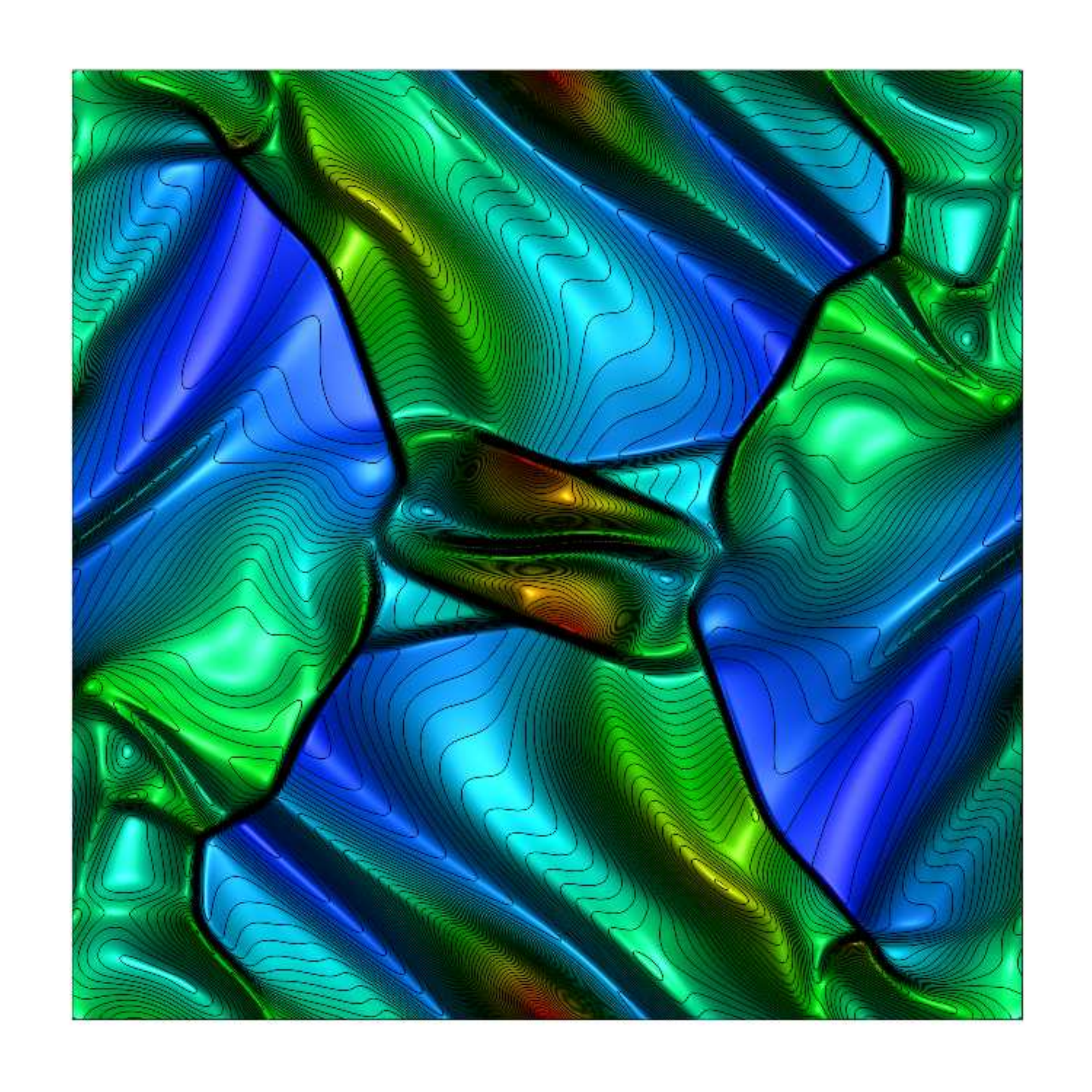} & 
\includegraphics[width=0.45\textwidth]{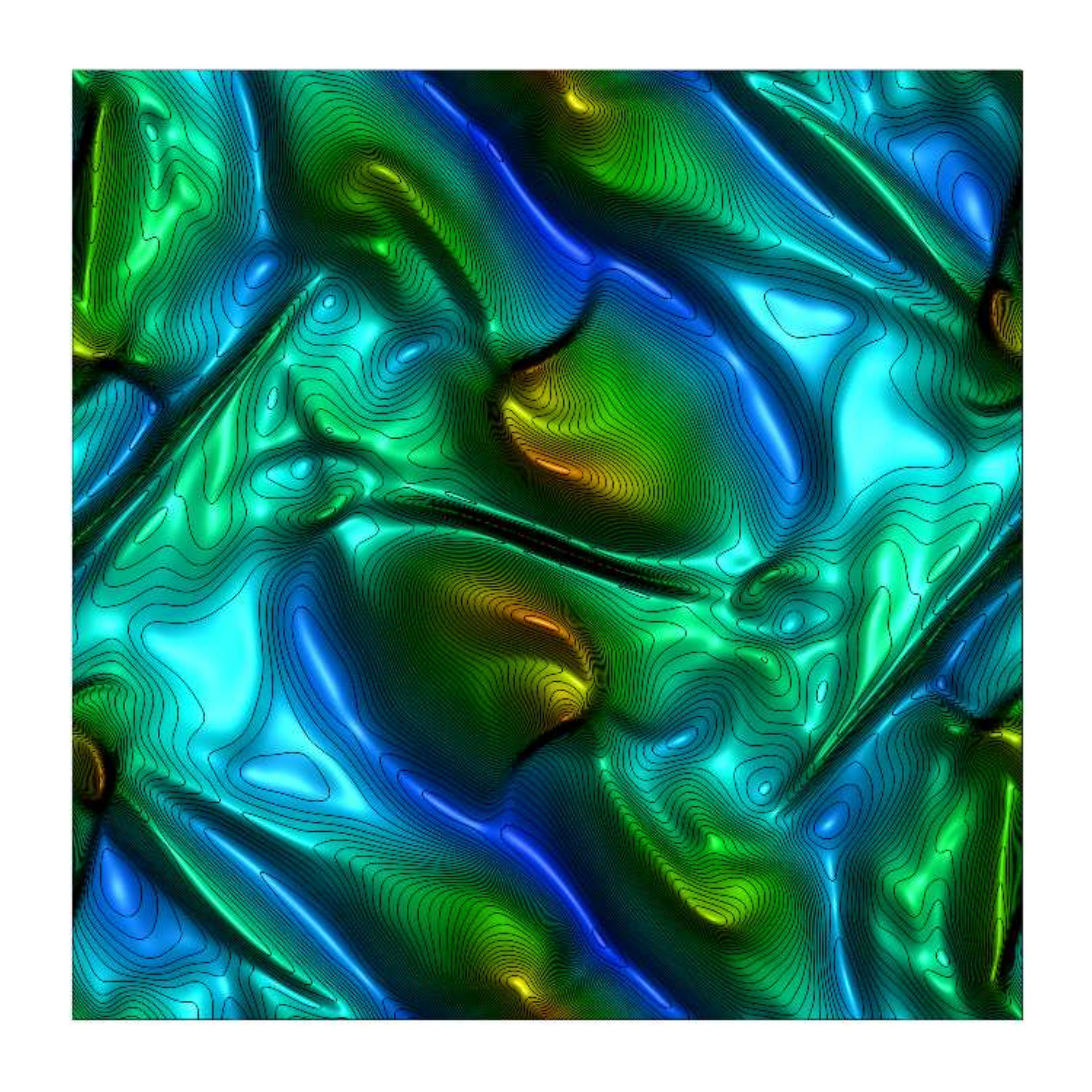}    
\end{tabular} 
\caption{Numerical solution obtained with the divergence-free semi-implicit finite volume method for the inviscid Orszag-Tang vortex system at time $t=0.5$ (top left), 
$t=2.0$ (top right), $t=3.0$ (bottom left) and $t=5.0$ (bottom right). 
56 equidistant contour lines of the pressure are shown in the interval $[0.5,6]$ are shown. }  
\label{fig.ot}
\end{center}
\end{figure}

Here we consider the very well-known Orszag-Tang vortex system for the ideal MHD equations, see \cite{OrszagTang,DahlburgPicone,PiconeDahlburg} for a detailed discussion 
of the underlying flow physics. The computational setup is the one used in \cite{JiangWu} and \cite{Dumbser2008} and is briefly summarized below.  
The computational domain under consideration is $\Omega = [0,2\pi]^2$ with four periodic boundary conditions. The initial conditions are given by 
$\rho = \gamma^2$, $\mathbf{v} = (-\sin(y),\sin(x),0)$, $p=\gamma$ and $\mathbf{B}=\sqrt{4 \pi} (-\sin(y),\sin(2x),0)$ with $\gamma = 5/3$. 
The computational domain is discretized with a uniform Cartesian mesh composed of $1000 \times 1000$ elements. The numerical results obtained with the SIFV scheme are
shown in Figure \ref{fig.ot} at times $t=0.5$, $t=2.0$, $t=3.0$ and $t=5.0$ and agree qualitatively well with those presented elsewhere in the literature, see e.g. 
\cite{Dumbser2008,balsarahlle2d,AMR3DCL,ADERdivB,HPRmodelMHD}. Also for this test case the average cost per element and time step was $3.0 \mu$s for the SIFV method.  
For comparison, the explicit second order accurate divergence-free Godunov-type scheme needed $2.4 \mu$s per element and time step, i.e. the average computational cost per
element and time step of the semi-implicit scheme is only about 25\% higher than for an analogous explicit method. Considering the fact that the semi-implicit scheme needs 
to solve $r_{\max}$ linear systems for the pressure in each time step (with $r_{\max}=2$ being the number of chosen Picard iterations), this means that the overhead due to the implicit 
discretization of the pressure is only very small for this test problem. In our view this is quite a remarkable result.  

\subsection{VRMHD current sheet and shear layer at low Mach number}
\label{sec.slcs}
\begin{figure}[!htbp]
\begin{center}
\begin{tabular}{cc} 
\includegraphics[width=0.45\textwidth]{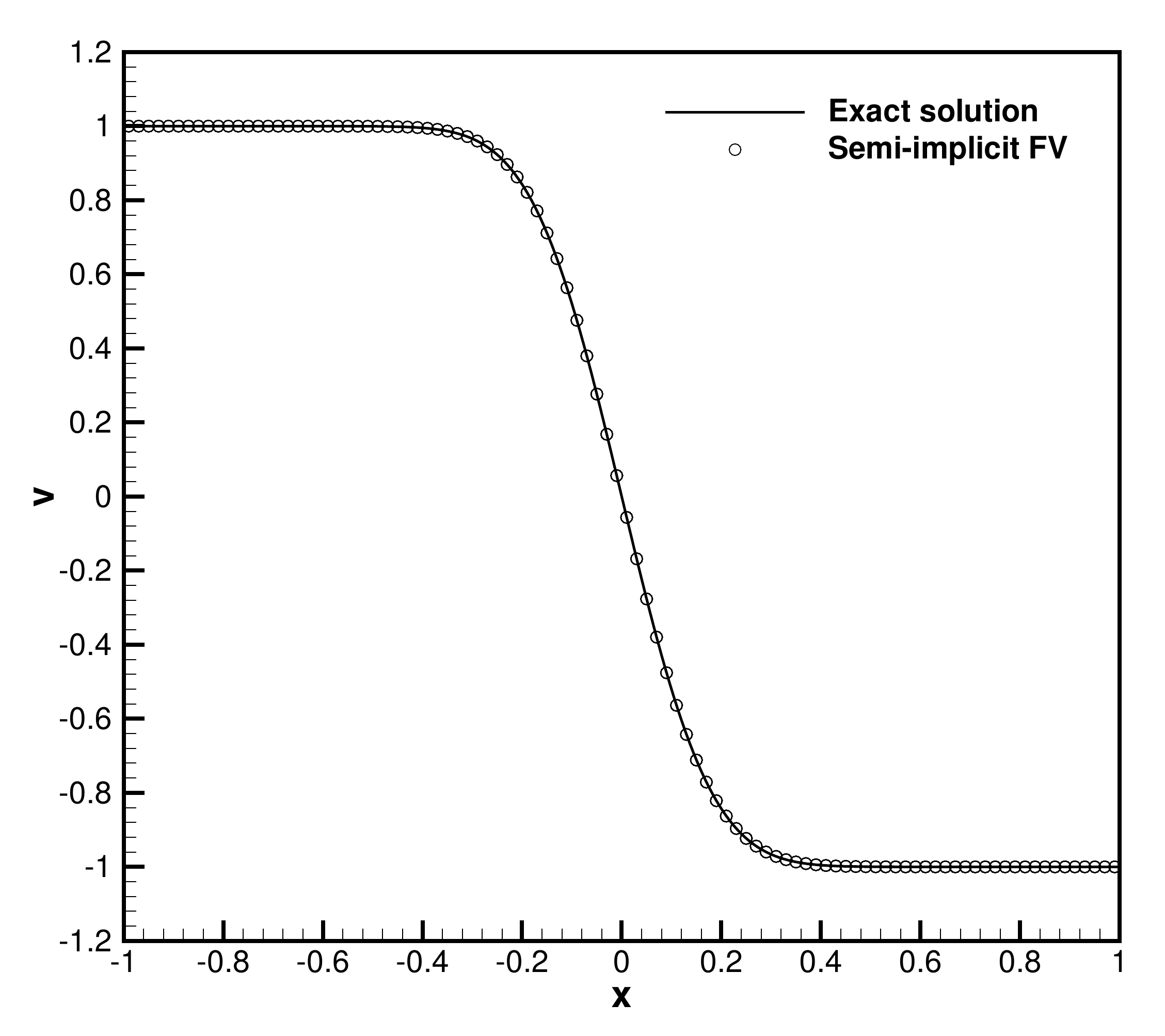} & 
\includegraphics[width=0.45\textwidth]{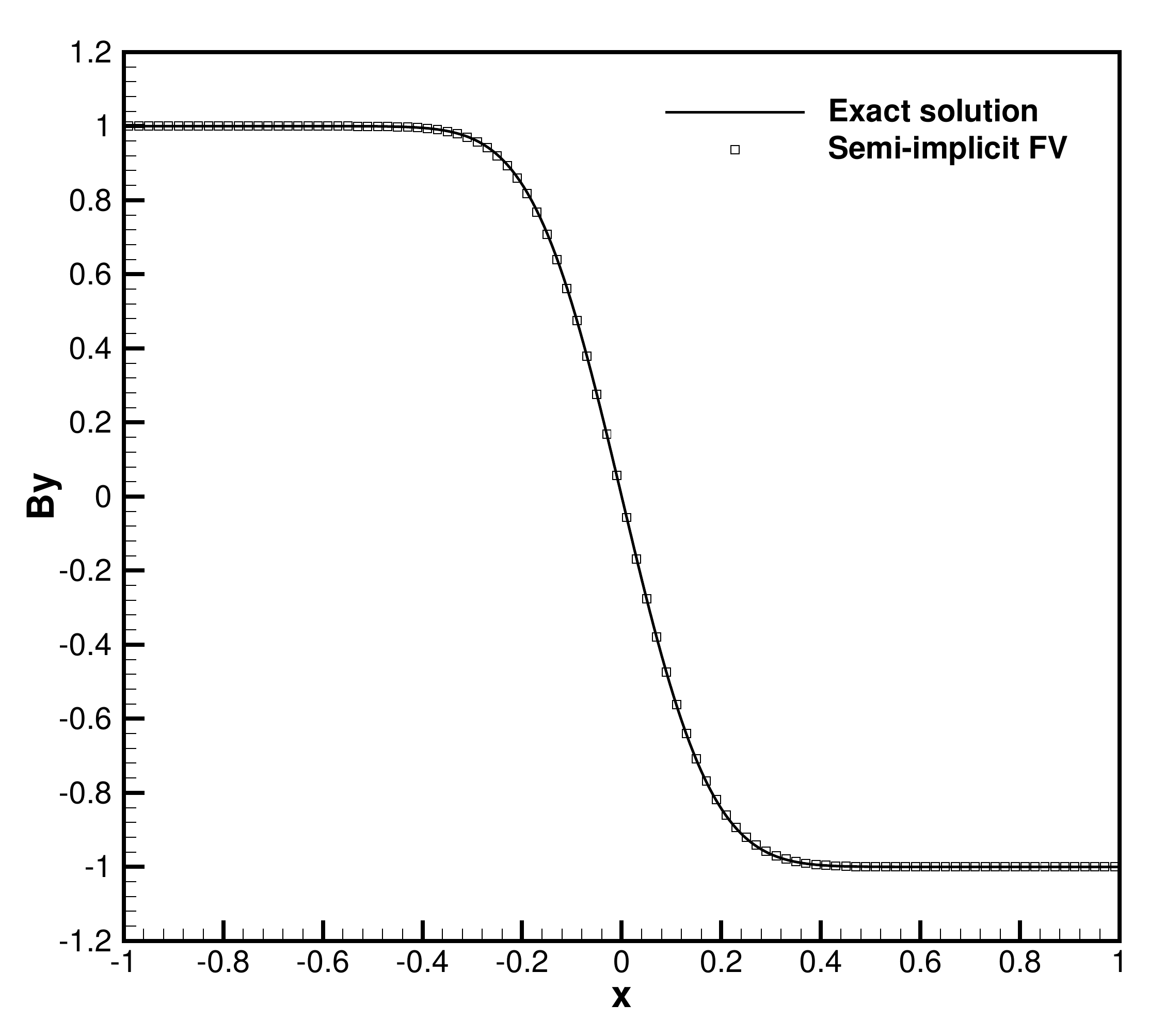}    
\end{tabular} 
\caption{Exact and numerical solution for the low Mach number shear layer (left) and the current sheet (right) at time $t=0.1$ solving the VRMHD equations with $\eta=\mu=10^{-1}$. } 
\label{fig.slcs}
\end{center}
\end{figure}

The current sheet and the simple shear layer (first problem of Stokes) are two very elementary test problems for the VRMHD equations, see e.g. \cite{Komissarov2007,ADERNSE}. 
Since our new semi-implicit finite volume scheme is particularly well-suited for low Mach number flows, we use the following initial conditions. In both cases, the density
and the fluid pressure are set to $\rho=1$ and $p=10^5$, respectively. For the shear layer, the initial magnetic field is zero and the velocity assumes the value 
$\mathbf{v}_L = (0,+1,0)$ for $x\leq 0$ and $\mathbf{v}_R=(0,-1,0)$ for $x>0$. The exact solution is given by (see \cite{boundaryLayer}): 
\begin{equation}
\label{eqn.stokes} 
  v(x,t) = - \textnormal{erf}\left( \halb \frac{x}{\sqrt{\mu t}} \right),  
\end{equation} 
We emphasize that this setup would be very challenging for an explicit solver due to the large value of the pressure and the resulting low Mach number. 
For the current sheet, the velocity is initialized with zero, while the magnetic field is
$\mathbf{B}_L = (0,+1,0)$ for $x\leq 0$ and $\mathbf{B}_R=(0,-1,0)$. The exact solution for $B_y$ is the same as the one given in \eqref{eqn.stokes} for the shear layer. 
In both cases the fluid parameters are $\eta=\mu=0.1$, $Pr=1$ and $c_v=1$. All simulations have been carried out until $t=0.1$ on the two-dimensional domain $\Omega=[-1,+1] \times [-0.1,+0.1]$ 
with periodic boundary conditions in $y$ direction and using a uniform Cartesian mesh of $100 \times 10$ elements. For this test we have deliberately chosen a 2D domain in order to check
our particular divergence-free implementation of the resistivity term at the aid of a discrete double curl. In Fig. \ref{fig.slcs} a scatter plot of the computational results obtained 
with the new divergence-free semi-implicit FV scheme is compared with the exact solution, where an excellent agreement can be observed for both cases. The scatter plot shows a clean 
one-dimensional behaviour, i.e. no spurious two-dimensional modes are introduced by the double curl operator. 

\subsection{VRMHD Orszag-Tang vortex}
\label{sec.vot}

\begin{figure}[!htbp]
\begin{center}
\begin{tabular}{cc} 
\includegraphics[width=0.45\textwidth]{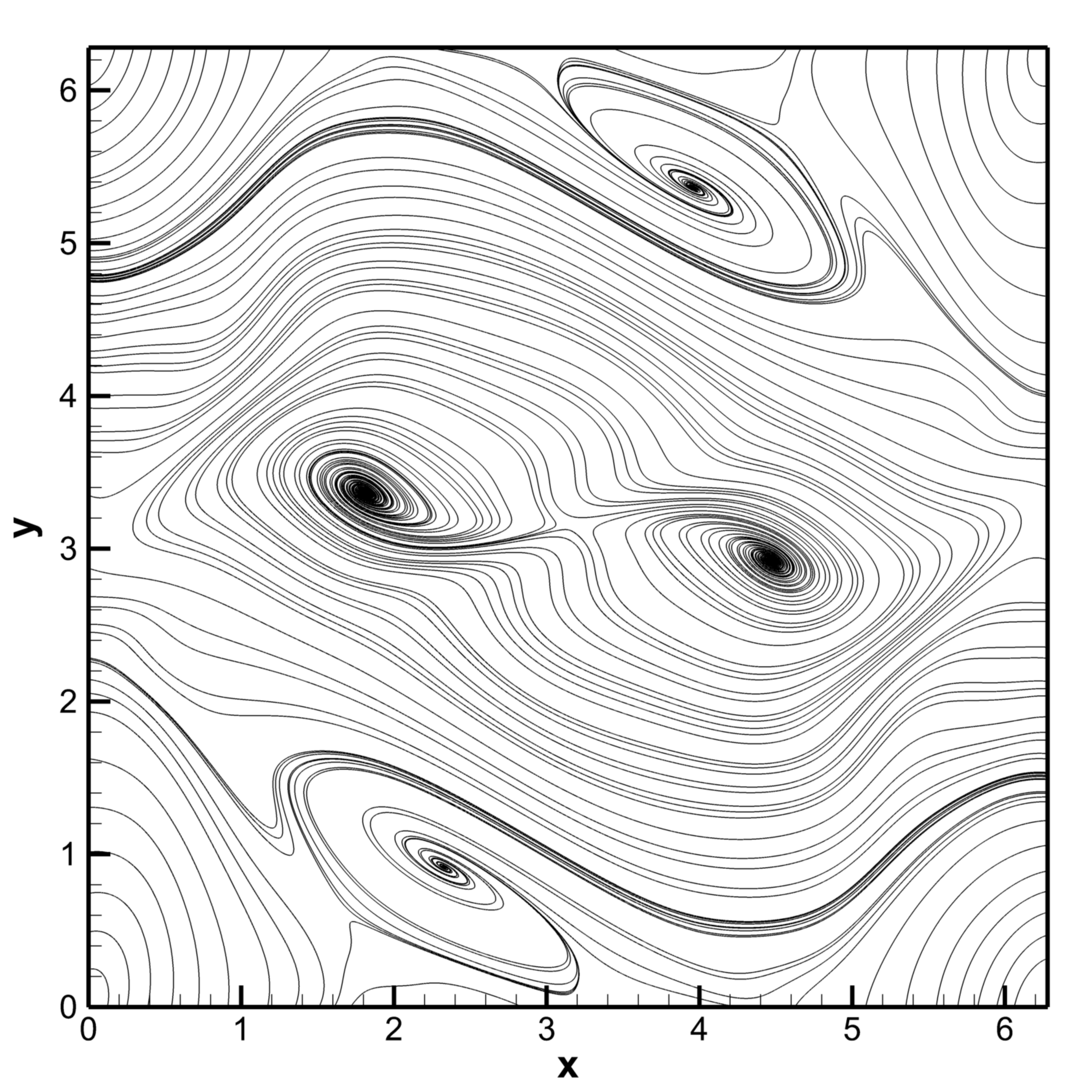}     & 
\includegraphics[width=0.45\textwidth]{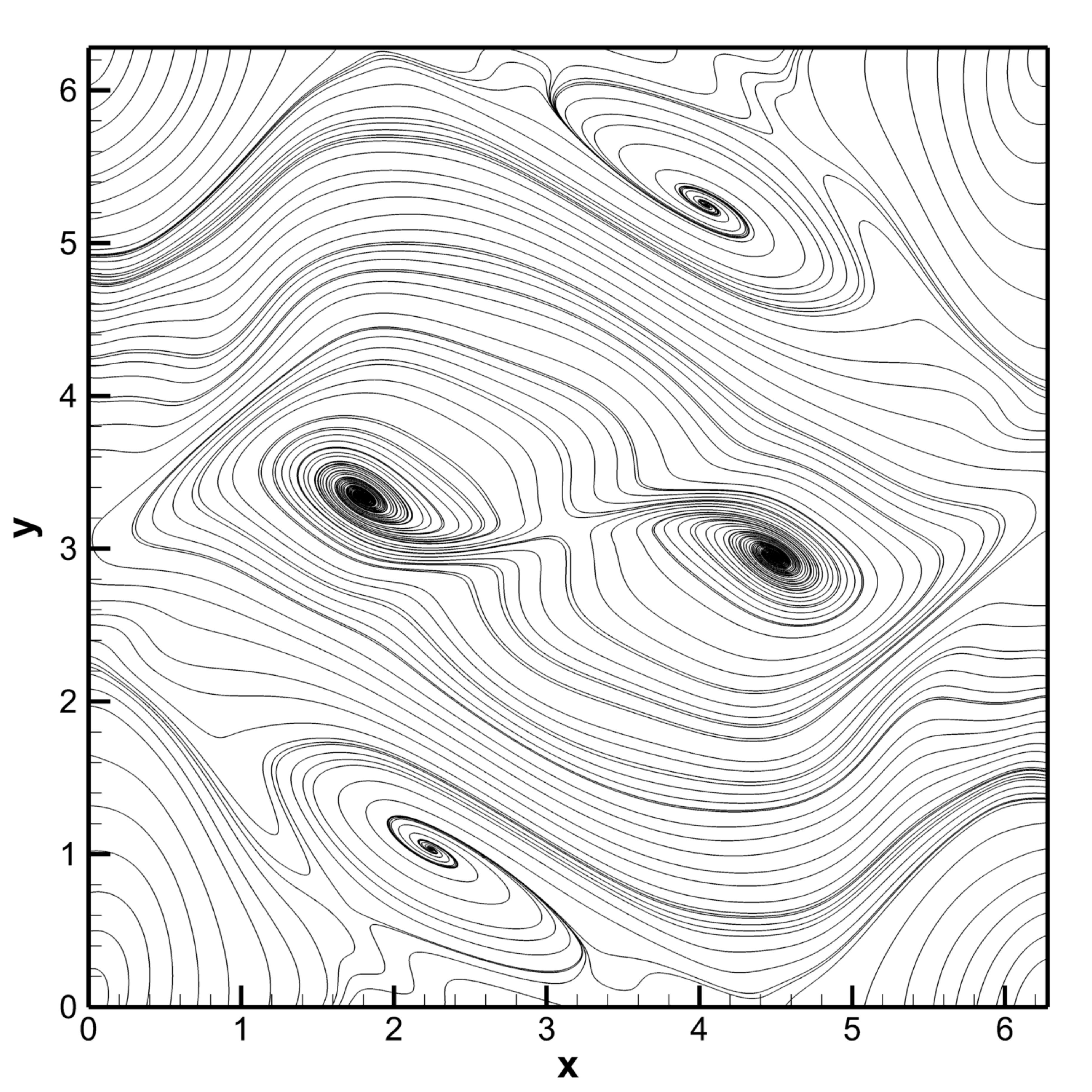}    \\   
\includegraphics[width=0.45\textwidth]{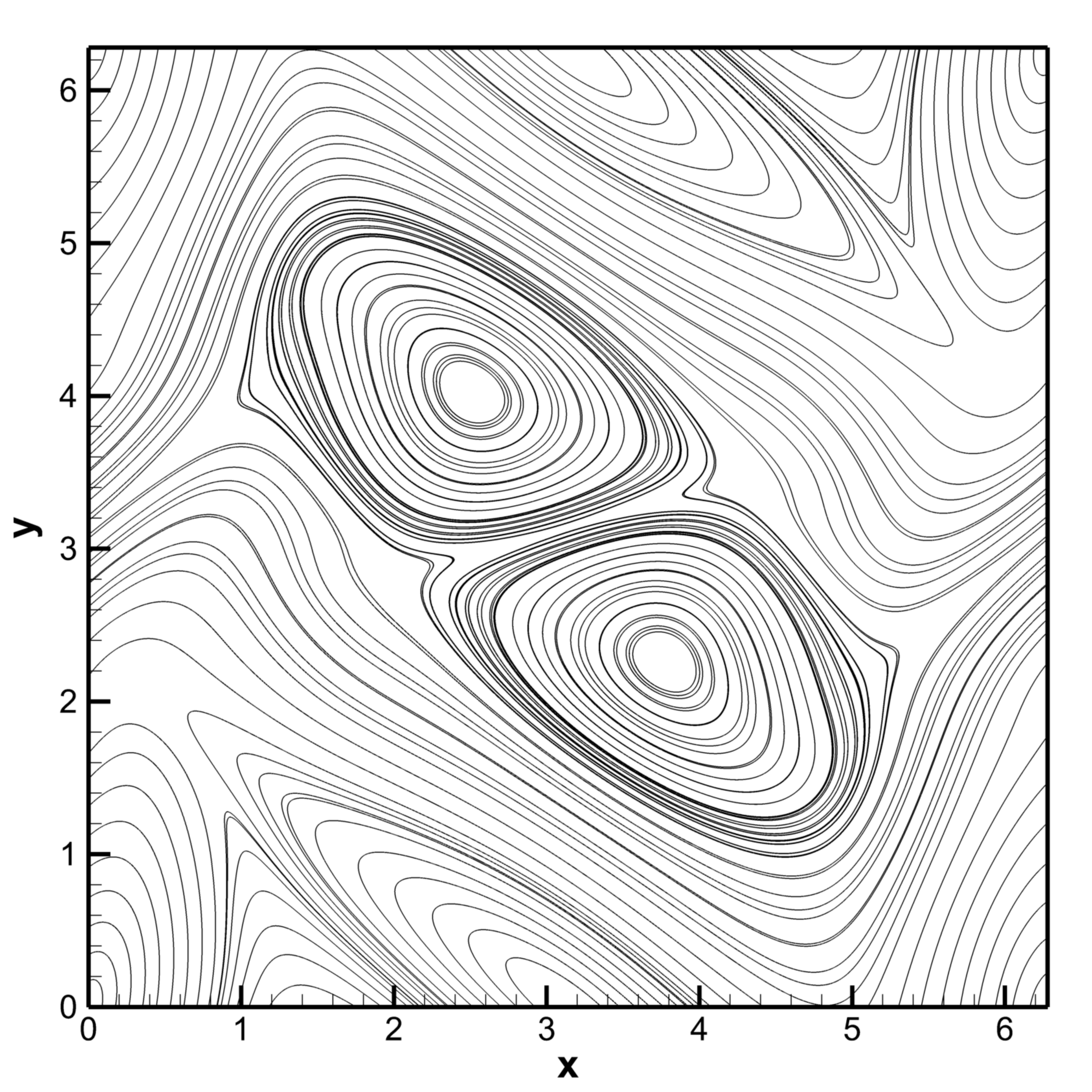}     & 
\includegraphics[width=0.45\textwidth]{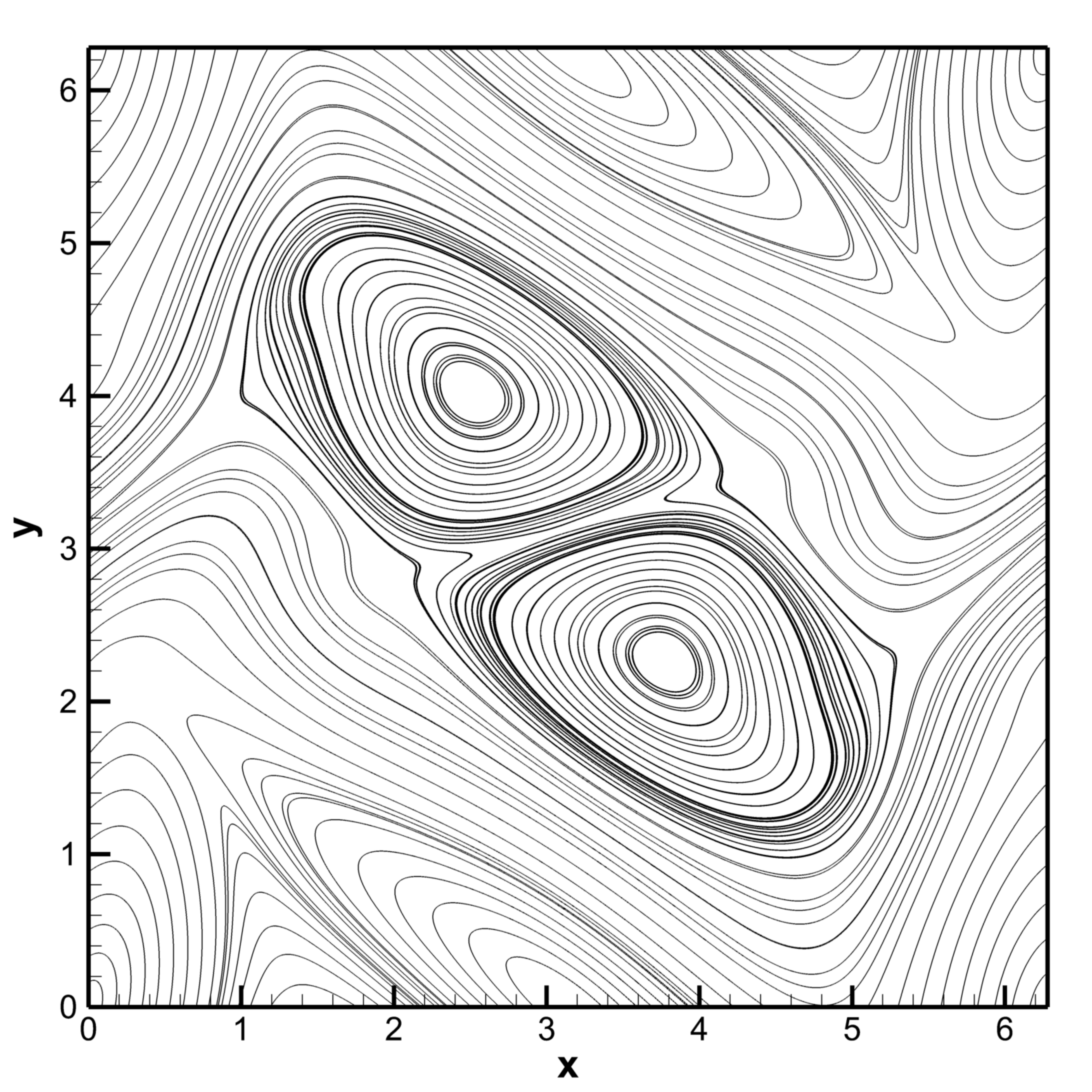}    \\   
\end{tabular} 
\caption{Reference solution (right) and numerical solution obtained with the new divergence-free semi-implicit finite volume method (left) for the viscous and resistive Orszag-Tang vortex 
($\eta=\mu=10^{-2}, Pr=1$) at time $t=2$. Velocity streamlines (top) and magnetic field lines (bottom). }
\label{fig.vot}
\end{center}
\end{figure}

In this subsection we solve the Orszag-Tang vortex system again, but this time using the viscous and resistive MHD equations (VRMHD). 
The fluid parameters are chosen as follows: $\gamma = \frac{5}{3}$, $\mu = \eta = 10^{-2}$, $c_v=1$ and a Prandtl number of $Pr=1$.  
The computational setup of this test case has been taken from \cite{WarburtonVRMHD} and \cite{ADERVRMHD} and is briefly summarized
below. The computational domain is again $\Omega=[0,2\pi]^2$ with four periodic boundary conditions, 
as in the inviscid case. The initial condition is given by $\rho=1$, $\mathbf{v}=\sqrt{4 \pi}(-\sin(y),\sin(x),0)$, 
$\mathbf{B}=(-\sin(y),\sin(2x), 0)$ and $p=\frac{15}{4} + \frac{1}{4} \cos(4x) + \frac{4}{5} \cos(2x) \cos(y) - \cos(x) \cos(y) + \frac{1}{4} \cos(2y)$. 
Simulations are carried out on a uniform Cartesian grid of $500 \times 500$ elements until a final time of $t=2$. The computational results obtained
with the SIFV scheme are shown in Fig. \ref{fig.vot}. They are also compared against a reference solution obtained in \cite{ADERVRMHD}
at the aid of a very high order accurate $P_NP_M$ scheme. Overall, we can note a good agreement between the two solutions. 
The average computational cost for this simulation was also $3\mu$s per element and time step, i.e. the scheme is able to update more than $3.33 \cdot 10^5$
zones per second on one single CPU core. 

\subsection{Kelvin-Helmholtz instability}
\label{sec.kh}

\begin{figure}[!htbp]
\begin{center}
\begin{tabular}{cc} 
\includegraphics[width=0.45\textwidth]{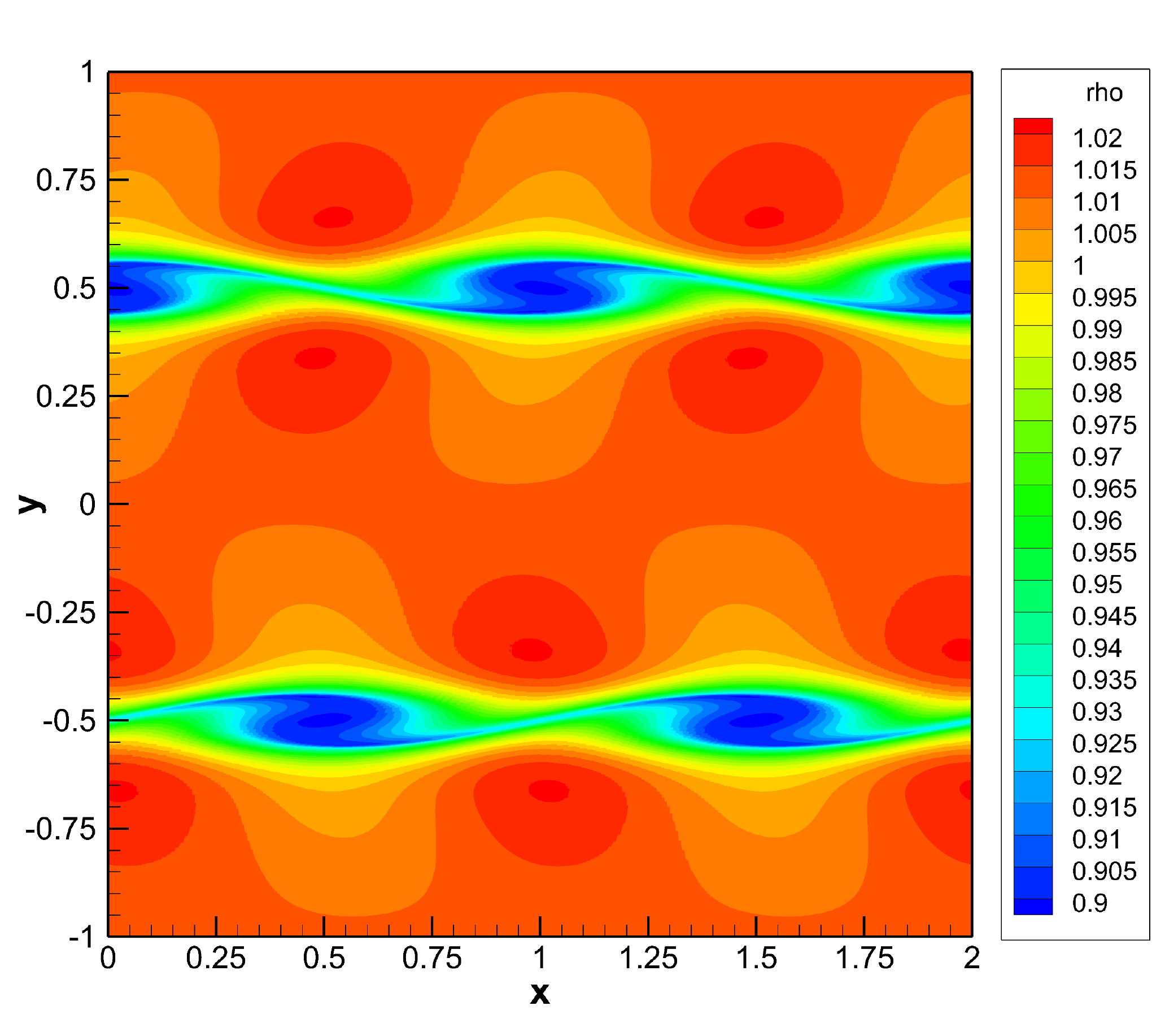}     & 
\includegraphics[width=0.45\textwidth]{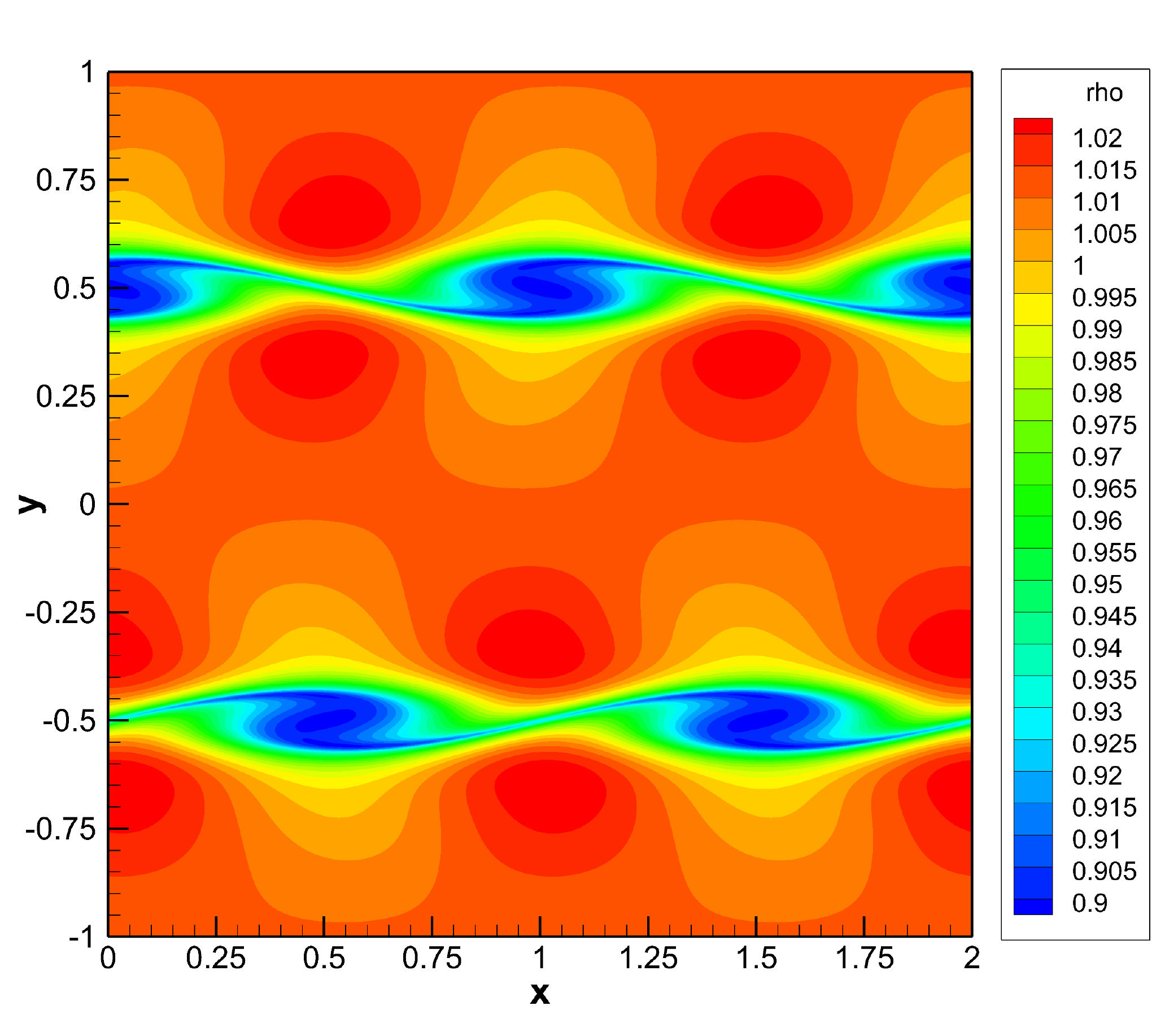} 
\end{tabular} 
\caption{Reference solution (right) and numerical solution obtained with the new divergence-free semi-implicit finite volume method (left) for the viscous and resistive MHD 
Kelvin-Helmholtz instability. The density contour levels are shown at the final time $t=4.0$.}
\label{fig.kh}
\end{center}
\end{figure}

In this test case we consider the same setup as presented in \cite{ADERVRMHD, HPRmodel} for the simulation of a Kelvin-Helmholtz instability developing in a viscous and resistive magnetized fluid. 
The initial condition is given by: $\rho=1$, $p=\frac{3}{5}$,  
\begin{equation*}
\mathbf{v}=\left(-\frac{1}{2}U_0 \tanh{\left(\frac{|y|-0.5}{a}\right)},\delta v \sin{(2 \pi x)}\sin{(\pi |y|)},0\right),
\end{equation*}
\begin{equation*}
	 \mathbf{B} = \left\{ 
	 \begin{array}{ccc}  
	 (B_0, 0, 0), & \textnormal{ if } & \halb + a < |y| < 1, \\ 
	 (B_0 \sin(\chi), 0, B_0 \cos(\chi)), & \textnormal{ if } & \halb - a < |y| < \halb + a, \\ 
	 (0, 0, B_0), & \textnormal{ if } & 0 < |y| < \halb - a,
	 \end{array}    	 
	   \right. 
\end{equation*}
with $ \chi = \frac{\pi}{2} \frac{y-0.5+a}{2a}$, $a=\frac{1}{25}$, $U_0=1$, $\delta v = 0.01$ and $B_0 = 0.07$. Furthermore $\gamma=\frac{5}{3}$, $\mu=\eta=10^{-3}$ and we neglect the heat conduction by setting $\lambda=0$. 
The computational domain is $\Omega=[0,2] \times [-1,1]$ using four periodic boundaries in all directions. For this test we use $1000 \times 1000$ elements and run the simulation up to $t=4s$. 
Figure \ref{fig.kh} shows the comparison between the numerical solution obtained with the proposed SIFV method and the one obtained in \cite{HPRmodel} using a high order explicit discontinuous  
Galerkin (DG) scheme for the solution of the VRMHD equations. A very good agreement can be observed also in this case that involves viscous and the resistive effects. The average computational 
cost for the new SIFV scheme was about $2.5 \mu$s per zone update. 

\section{Conclusions}
\label{sec.concl} 
In this paper we have presented a new divergence-free semi-implicit finite volume method for the simulation of the ideal, viscous and resistive magnetohydrodynamics equations with 
general equation of state. 

The split discretization of the nonlinear convective and viscous terms on the main grid combined with our very particular discretization of the pressure subsystem on a staggered grid allows us 
to reduce the final problem to the solution of a mildly nonlinear system for the fluid pressure, which can be efficiently solved by the (nested) Newton-type technique of Casulli et al. 
\cite{BrugnanoCasulli,BrugnanoCasulli2,CasulliZanolli2010,CasulliZanolli2012}. The linear part of the mildly nonlinear system is given by a symmetric and positive semi-definite M-matrix, 
which is a very remarkable property for a semi-implicit time discretization of the MHD equations. 
The nonlinearity in our mildly nonlinear system resides only on the diagonal and is contained in the equation of state that needs to provide the specific energy $e=e(p,\rho)$ as a function 
of the fluid pressure and the density. The EOS must must be a non-negative non-decreasing function of $p$ (for a given density) and its partial derivative w.r.t. $p$ must be a function of 
bounded variation. 
For linear equations of state like the ideal gas EOS, the entire pressure system becomes linear and can therefore be solved in one single Newton iteration. 
The unknown kinetic energy at the new time level as well as the specific enthalpies are updated easily with a simple Picard process, following the suggestion of \cite{CasulliZanolli2010}. 
Once the pressure is known at the new time level, the momentum and total energy density can be readily obtained via a conservative update formula. 

The magnetic field in our new SIFV scheme is also discretized on the staggered mesh, following the ideas of Balsara et al. \cite{balsarahlle2d,balsarahllc2d,balsarahlle3d,BalsaraMultiDRS,MUSIC1,MUSIC2} 
on exactly divergence-free schemes for MHD and multi-dimensional Riemann solvers. In our method the resistive terms in the induction equation are discretized using a discrete double curl formulation,  which assures that the scheme remains exactly divergence free also in the non-ideal (resistive) case.

The time step of our new method is only restricted by the fluid velocity and the speed of the Alfv\'en waves, but not by the speed of sound. 
Therefore, our scheme is particularly well-suited for \textbf{low Mach number flows}. For example, in the low Mach number magnetic field loop advection test presented in Section \ref{sec.fieldloop}, 
our new semi-implicit method was more than \textbf{50 times faster} compared to a comparable explicit divergence-free second-order accurate Godunov-type finite volume scheme. 
Nevertheless, extensive numerical experiments have shown that our new pressure-based solver performs very well also for high Mach number flows with shock waves and other flow discontinuities. 
We have also compared the computational cost of the new SIFV scheme 
with the cost of a standard second-order Godunov-type scheme for MHD using the same code basis and the same computer, in order to get a fair comparison. For example, for the Orszag-Tang 
vortex problem shown in Section \ref{sec.ot} the average cost per element and time step of the explicit scheme was about $2.4 \mu$s, while it was about $3.0\mu$s for the semi-implicit method, 
i.e. despite the necessary solution of a linear system for the pressure in each of the two Picard iterations of the SIFV scheme, the semi-implicit method was only 25\% more expensive 
than a fully explicit discretization. This means that we have a \textit{very low overhead} due to the implicit discretization of the pressure subsystem, which in our opinion is also a 
remarkable result. 

Future work will consist in an extension of the present approach to general unstructured meshes in multiple space dimensions and to higher order of accuracy at the aid of staggered 
semi-implicit discontinuous Galerkin (DG) finite element schemes, following the ideas outlined in \cite{DumbserCasulli,TavelliDumbser2016,TavelliDumbser2017,FambriDumbser,AMRDGSI}. 
In the near future we also plan an extension of this new family of efficient semi-implicit finite volume schemes to the unified Godunov-Peshkov-Romenski (GPR) model of continuum mechanics 
\cite{PeshRom2014,HPRmodel,HPRmodelMHD} and to the Baer-Nunziato model of compressible multi-phase flows \cite{BaerNunziato1986,SaurelAbgrall,SaurelAbgrall2}, where low Mach number 
problems are particularly important due to the simultaneous presence of two different phases.

\section*{Acknowledgements}

The authors would like to thank S.A.E.G. Falle for providing the exact Riemann solver for the ideal MHD equations. 
The research presented in this paper was partially funded by the European Research Council (ERC) under the European Union's Seventh Framework Programme (FP7/2007-2013) within the research project 
\textit{STiMulUs}, ERC Grant agreement no. 278267 and by the European Union's Horizon 2020 Research and Innovation Programme under the project \textit{ExaHyPE}, grant no. 671698 (call FETHPC-1-2014).  

\bibliography{SIMHD}
\bibliographystyle{plain}

\end{document}